\documentclass[a4paper,12pt]{article}
\usepackage{dcolumn}
\usepackage{bm}
\usepackage[dvips]{graphicx}
\usepackage{amsmath}
\usepackage{hyperref}
\usepackage{amssymb}
\usepackage{amsfonts}
\usepackage{ucs}
\usepackage[utf8x]{inputenc}
\usepackage[T1]{fontenc}
\usepackage{url}
\usepackage{psfrag}
\usepackage{dsfont}
\usepackage{color}

\newtheorem{proposition}{\textbf{Proposition}}

\newtheorem{corollary}{\textbf{Corollary}}
\newtheorem{remark}{\textbf{Remark}}
\newtheorem{theorem}{\textbf{Theorem}}
\newtheorem{lemma}{\textbf{Lemma}}
\newtheorem{definition}{\textbf{Definition}}

\newcommand{\noautospacebeforefdp}{}
\newcommand{\autospacebeforefdp}{}
\newcommand{\NoAutoSpaceBeforeFDP}{}
\newcommand{\AutoSpaceBeforeFDP}{}

\def\dx{{\rm d}x}

\def\del  {\partial}
\def\eps{\varepsilon}

\def\R{\mathbb{R}}
\def\r{\mathbb{R}}

\def\N{\mathbb{N}}
\def\C{\mathbb{C}}

\def\dx{{\rm d}x}
\def\dy{{\rm d}y}

\def\dm{{\rm d}m}

\begin{document}
\title{Dirichlet-to-Neumann or Poincaré-Steklov operator on fractals described by $d$-sets}

\author{KEVIN ARFI\footnote{kevin.arfi@student.ecp.fr},  ANNA ROZANOVA-PIERRAT\footnote{Laboratoire Math\'ematiques et Informatique Pour la Complexit\'e et les Syst\`emes,
Centrale Sup\'elec, Universit\'e Paris-Saclay, Grande Voie des Vignes,
Ch\^atenay-Malabry, France, 
anna.rozanova-pierrat@centralesupelec.fr}}

\maketitle

\begin{abstract}
In the framework of the Laplacian transport, described by a Robin boundary value problem in an exterior domain in $\R^n$, we generalize the definition of the Poincar\'e-Steklov operator to $d$-set boundaries, $n-2< d<n$, and give its spectral properties to compare to  the spectra of the interior domain and also of a truncated domain, considered as an approximation of the exterior case. 
The well-posedness of the Robin boundary value problems for the truncated and exterior domains is given in the general framework of $n$-sets. The results are obtained thanks to a generalization of the continuity and compactness properties of the trace and extension operators in Sobolev, Lebesgue and Besov spaces, in particular, by  a generalization of the  classical Rellich-Kondrachov Theorem of compact embeddings for $n$ and $d$-sets.
\end{abstract}

\section{Introduction}\label{SecIntro}

Laplacian transports to and across irregular and fractal interfaces are ubiquitous in nature and industry: properties of rough electrodes in electrochemistry,
heterogeneous catalysis, steady-state transfer across biological membranes (see~\cite{GREBENKOV-2006,GREBENKOV-2007,FILOCHE-2008,GREBENKOV-2004}
and references therein). To model it there is a usual interest to consider truncated domains as an approximation of the exterior unbounded domain case. 

\begin{figure}[!ht]
 \begin{center}
    \includegraphics[width=7in]{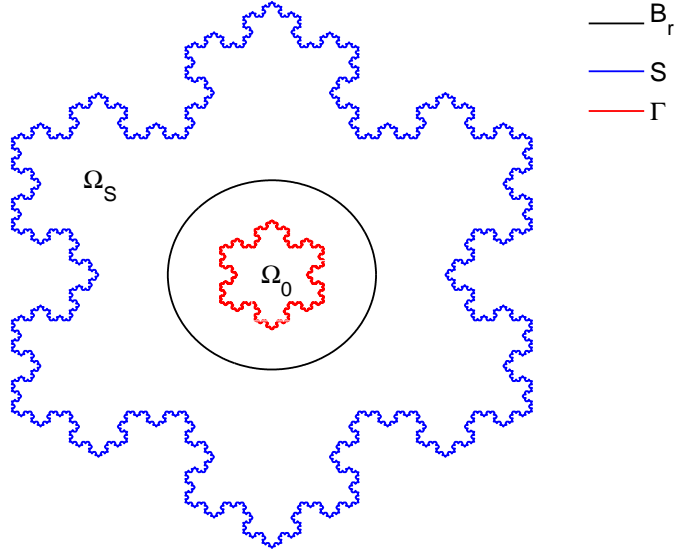}
    \end{center} 
    \caption{\label{FigSchema} Example of the considered domains: $\Omega_0$ (the Von Koch snowflake) is the bounded domain, bounded by a compact boundary $\Gamma$, which is a $d$-set (see Definition~\ref{Defdset}) with $d=\log 4/ \log 3>n-1=1$. The truncated domain $\Omega_S$ is between the boundary $\Gamma$ and the boundary $S$ (presented by the same Von Koch fractal as $\Gamma$). The boundaries $\Gamma$ and $S$ have no an intersection and here are separated by the boundary of a ball $B_r$ of a radius $r>0$. The domain, bounded by $S$, is called $\Omega_1= \overline{\Omega}_0\cup \Omega_S$, and the exterior domain is $\Omega=\R^n\setminus \overline{\Omega}_0$.}
  \end{figure}

Let $\Omega_0$ and $\Omega_1$ be two bounded domains in $\R^n$  with disjoint boundaries 
 $\del \Omega_0\cap \del \Omega_1=\varnothing$, denoted by $\Gamma$ and $S$ respectively, such that $\overline{\Omega}_0\subset \Omega_1$.
 Thus, in this paper, we consider two types of domains constructed on $\Omega_0$:
\begin{enumerate}
 \item the unbounded exterior domain to $\Omega_0$, denoted by $\Omega=\R^n\setminus \overline{\Omega}_0$;
 \item a bounded, truncated by a boundary $S$, truncated  domain  $\Omega_S=(\R^n\setminus \overline{\Omega}_0)\cap \Omega_1$.
\end{enumerate}
 Let us notice that $\Gamma\cup S=\del \Omega_S$ (for the unbounded case $S=\varnothing$ and $\del \Omega=\Gamma$), see Fig.~\ref{FigSchema}. As $\Omega_0$ is bounded, its boundary $\Gamma$ is supposed compact. 
The phenomenon of Laplacian transport to $\Gamma$ can be described by the following boundary value problem:
\begin{equation}\label{Problem}
 \begin{split}
&  -\Delta u = 0, \quad x\in \Omega_S\hbox{ or } \Omega,\\
&\lambda u + \partial_\nu u = \psi \quad \mbox{on } \Gamma,\\
& u = 0 \quad  \mbox{on } S,
 \end{split}
\end{equation}
where $\partial_\nu u$ denotes the normal derivative of $u$, in some appropriate sense,  $\lambda \in [0, \infty)$ is the resistivity of the boundary and $\psi \in L_2(\Gamma)$.
For $S=\varnothing$ we impose Dirichlet boundary conditions at infinity. The case of a truncated domain $\Omega_S$ corresponds to an approximation of the exterior problem  in the sense of Theorem~\ref{theorem/convOpRobin}. 

When $\del \Omega$ is regular ($C^\infty$ or at least Lipschitz), it is well-known~\cite{LIONS-1972,MARSCHALL-1987} how to define the trace of $u\in H^1(\Omega)$ and the normal derivative  $\partial_\nu u$ on $\Gamma$. The properties of the Poincar\'e-Steklov  or the Dirichlet-to-Neumann operator, defined at manifolds with $C^\infty$-boundaries are also well-known~\cite{GIROUARD-2014,TAYLOR-1996}. In the aim to generalize the Poincar\'e-Steklov operator to $d$-sets with $n-2< d <n$ (the case $n-1<d<n$ contains the self-similar fractals), we firstly study the most general context (see Section~\ref{SecWellPosed}), 
   when the problem~(\ref{Problem}) is well-defined and its bounded variant (physically corresponding to a source at finite distance) can be viewed as an approximation of the unbounded case (corresponding to a source at infinity). The main extension and trace theorems, recently obtained in the framework of $d$-sets theory, are presented and discussed in Section~\ref{SecFuncAnT}. They allow us to generalize the known  properties of the trace and extension operators on the $(\eps,\delta)$-domains~\cite{JONES-1981,WALLIN-1991} (see Theorem~\ref{ThContTrace}) to a more general class of $n$-sets, called  admissible domains (see Definition~\ref{DefAdmis}), and update for admissible domains the classical Rellich-Kondrachov theorem (see Theorems~\ref{ThCSEnSET} and~\ref{ThAnnaEmbF} for $d$-sets). Actually, we state that the compactness of a Sobolev embedding to a Sobolev space does not depend on the boundness of the domain, but it is crucial for the embeddings in the Lebesgue spaces. Hence, a trace operator $H^1(\Omega)\to L_2(\del \Omega)$ mapping the functions defined  on a domain $\Omega$ to their values on the boundary $\del \Omega$ (or on any part $D$ of $\Omega$, $H^1(\Omega)\to L_2(D)$) is compact if and only if the boundary $\del \Omega$ (or the part $\overline{D}$) is compact.  
   
 After a short survey in Section~\ref{SecDtoNReg} of known results on the spectral properties of the Poincar\'e-Steklov operator for a bounded domain, 
 we introduce the Poincar\'e-Steklov operator $A$ on a compact $d$-set boundary $\Gamma$ of an admissible bounded domain $\Omega_0$.
 Since  $\Gamma$ (see Fig.~\ref{FigSchema}) can be viewed not only as the boundary of $\Omega_0$, but also as the boundary of the exterior domain $\Omega$ and of its truncated domain $\Omega_S$, we also introduce the Poincar\'e-Steklov operator $A$ on $\Gamma$ for the exterior and trucated cases and relate their spectral properties (see Section~\ref{SecDtoNFrac}).
  In all cases, the Poincar\'e-Steklov operator $A$ can be defined as a positive self-adjoint operator on $L_2(\Gamma)$, and $A$ has a discrete spectrum if and only if the boundary $\Gamma$ is compact.   The two dimensional case differs from the case of $\R^n$ with $n\ge 3$ by the functional reason (see Subsection~\ref{SubsFuncSpW}) and gives different properties of the point spectrum of $A$ (see Theorem~\ref{ThSPECTR}). In particular, in the exterior case $A$ for $n=2$ and $n\ge 3$ has different domains of definition (see Proposition~\ref{Prop1f} in Section~\ref{SecFinal}). 
  
   Specially, for the case of a $d$-set $\Gamma$ (see Theorems~\ref{ThDtoNInt},~\ref{PropDtNn2} and~\ref{ThDtNcompResolventeExt}), we justify the method, developed in~\cite{GREBENKOV-2006}, true for smooth boundaries, to find the total flux $\Phi$ across the interface $\Gamma$ using the spectral decomposition of $1_\Gamma$ (belonging to the domain of $A$ by Proposition~\ref{Prop1f}) on the basis of eigenfunctions of the Dirichlet-to-Neumann operator $(V_k)_{k\in \N}$ in $L_2(\Gamma)$ and its eigenvalues $(\mu_k)_{k\in \N}$:
\begin{equation}
\label{eq/Phi}
\Phi \propto \sum_k \frac{\mu_k (\mathds{1}_\Gamma , V_k)_{L_2(\Gamma)}^2}{1 + \frac{\mu_k}{\lambda}}.
\end{equation}

 \section{Continuity and compactness of the extension and trace operators on $d$-sets}\label{SecFuncAnT}
          
    Before to proceed to the generalization results, let us define the main notions and  explain the functional context of $d$-sets. For instance, for the well-posedness result of problem~(\ref{Problem}) on ``the most general'' domains $\Omega$ in $\R^n$, we need to be able to say that for this $\Omega$ 
  the extension operator $\rm{E}:$ $H^1(\Omega)\to H^1(\R^n)$ is continuous and the trace operator (to be defined, see Definition~\ref{DefGTrace})
  $\mathrm{Tr}: H^1(\Omega)\to \rm{Im(Tr}(H^1(\Omega)))\subset L_2(\del \Omega)$ is continuous and surjective.
          
Therefore, let us introduce the existing results about traces and extension domains in the framework of Sobolev spaces.
\begin{definition}\textbf{($W_p^k$-extension domains)}
 A domain $\Omega\subset \R^n$ is called a $W_p^k$-extension domain ($k\in \N^*$) if there exists a bounded linear extension operator $E: W_p^k(\Omega) \to W^k_p(\R^n)$. This means that for all $u\in W_p^k(\Omega)$ there exists a $v=Eu\in  W^k_p(\R^n)$ with $v|_\Omega=u$ and it holds
 $$\|v\|_{W^k_p(\R^n)}\le C\|u\|_{W_p^k(\Omega)}\quad \hbox{with a constant } C>0.$$
\end{definition}

The classical results of Calderon-Stein~\cite{CALDERON-1961,STEIN-1970} say that every Lipschitz domain $\Omega$ is an extension domain for $W_p^k(\Omega)$ with $1\le p\le \infty$, $k\in \N^*$.

This result was generalized by Jones~\cite{JONES-1981} in the framework of $(\eps,\delta)$-domains:
\begin{definition}\label{DefEDD}\textbf{($(\eps,\delta)$-domain~\cite{JONES-1981,JONSSON-1984,WALLIN-1991})}
An open connected subset $\Omega$ of $\R^n$ is an $(\eps,\delta)$-domain, $\eps > 0$, $0 < \delta \leq \infty$, if whenever $x, y \in \Omega$ and $|x - y| < \delta$, there is a rectifiable arc $\gamma\subset \Omega$ with length $\ell(\gamma)$ joining $x$ to $y$ and satisfying
\begin{enumerate}
 \item $\ell(\gamma)\le \frac{|x-y|}{\eps}$ and
 \item $d(z,\del \Omega)\ge \eps |x-z|\frac{|y-z|}{|x-y|}$ for $z\in \gamma$. 
\end{enumerate}
\end{definition}


This kind of domains are also called locally uniform domains~\cite{HERRON-1991}. 
Actually, bounded locally uniform domains, or  bounded $(\eps,\delta)$-domains, are equivalent (see~\cite{HERRON-1991} point 3.4) to the uniform domains, firstly defined by Martio and Sarvas in~\cite{MARTIO-1979}, for which there are no more restriction $|x-y|<\delta$  (see Definition~\ref{DefEDD}).

Thanks to Jones~\cite{JONES-1981}, it is known that any $(\eps,\delta)$-domain in $\R^n$ is a $W_p^k$-extension domain for all $1\le p\le\infty$ and $k\in \N^*$. Moreover, for a bounded finitely connected domain $\Omega\subset \R^2$, Jones~\cite{JONES-1981} proved that 
$\Omega$ is a $W_p^k$-extension domain ($1\le p\le\infty$ and $k\in \N^*$) if and only if $\Omega$ is an $(\eps,\infty)$-domain for some $\eps>0$, if and only if the boundary $\del \Omega$ consists of finite number of points and quasi-circles.
However, it is no more true for $n\ge3$, $i.e.$ there are $W_p^1$-extension domains which are not locally uniform~\cite{JONES-1981} (in addition, an $(\eps,\delta)$-domain in $\R^n$ with $n\ge 3$ is not necessary a quasi-sphere).

To discuss general properties of  locally uniform domains, let us introduce Ahlfors $d$-regular sets or $d$-sets:
\begin{definition}\label{Defdset}\textbf{(Ahlfors $d$-regular set or $d$-set~\cite{JONSSON-1984,WALLIN-1991,JONSSON-1995})} 
Let $F$ be a Borel subset of $\R^n$  and $m_d$ be the $d$-dimensional Hausdorff measure,  $0<d\le n$.  The set $F$ is called a $d$-set, if there exist  positive constants $c_1$, $c_2>0$,
 \begin{equation*}
 c_1r^d\le m_d(F\cap B_r(x))\le c_2 r^d, \quad \hbox{ for  } ~ \forall~x\in F,\; 0<r\le 1,
 \end{equation*}
where $B_r(x)\subset \R^n$ denotes the Euclidean ball centered at $x$ and of radius~$r$.
\end{definition}
%
Henceforth, the boundary $\Gamma$ is a $d$-set endowed with the $d$-dimensional Hausdorff measure, and $L_p(\Gamma)$ is defined with respect to this measure as well.

From~\cite{WALLIN-1991}, it is known that
\begin{itemize}
 \item All $(\eps,\delta)$-domains in $\R^n$ are $n$-sets ($d$-set with $d=n$):
 $$\exists c>0\quad \forall x\in \overline{\Omega}, \; \forall r\in]0,\delta[\cap]0,1] \quad \mu(B_r(x)\cap \Omega)\ge C\mu(B_r(x))=cr^n,$$
 where $\mu(A)$ denotes the Lebesgue measure of a set $A$. This property is also called the measure density condition~\cite{HAJLASZ-2008}. Let us notice that an $n$-set 
$\Omega$ cannot be ``thin'' close to its boundary $\del \Omega$.
 \item If $\Omega$ is  an $(\eps,\delta)$-domain and $\del \Omega$ is a $d$-set ($d<n$) then $\overline{\Omega}=\Omega\cup \del \Omega$ is an $n$-set.
\end{itemize}
In particular, a Lipschitz domain $\Omega$ of $\R^n$ is an $(\eps,\delta)$-domain and also an $n$-set~\cite{WALLIN-1991}. But not every $n$-set is an $(\eps,\delta)$-domain: adding an in-going cusp to an $(\eps,\delta)$-domain we obtain an $n$-set which is not an $(\eps,\delta)$-domain anymore. 
Self-similar fractals (e.g., von Koch's snowflake domain) are examples of $(\eps,\infty)$-domains with the $d$-set boundary~\cite{CAPITANELLI-2010,WALLIN-1991}, $d>n-1$.
From~\cite{JONSSON-1984} p.39, it is also known that all closed $d$-sets with $d>n-1$ preserve Markov's local inequality:
\begin{definition}\textbf{(Markov's local inequality)}
A closed subset $V$ in $\R^n$ preserves Markov's local inequality if for every fixed  $k\in \N^*$, there exists a constant $c=c(V,n,k) > 0$, such that 
$$
\max_{V\cap \overline{B_r(x)}} |\nabla P | \le \frac{c}{r}\max_{V\cap \overline{B_r(x)}}|P|
$$
for all polynomials $P \in \mathcal{P}_k$ and all closed balls $\overline{B_r(x)}$, $x \in V$ and $0 < r \le 1$.
\end{definition}
%

For instance, self-similar sets that are not  subsets of any $(n-1)$-dimensional subspace of $\R^n$, the closure of a domain $\Omega$ with Lipschitz boundary and also $\R^n$ itself preserve Markov's local inequality (see Refs.~\cite{WALLIN-1991,JONSSON-1997}). 
The geometrical characterization of sets preserving Markov's local inequality was initially given in~\cite{JONSSON-1984-1} (see Theorem 1.3) and can be simply interpreted as sets which are not too flat anywhere. It can be  illustrated by the following theorem of Wingren~\cite{WINGREN-1988}:
\begin{theorem}
A closed subset $V$ in $\R^n$ preserves Markov's local inequality if and only if there exists a constant $c>0$ such that for every ball
$B_r(x)$ centered in $x\in V$ and with the radius $0 < r \le 1$, there are $n + 1$
affinely independent points $y_i \in V\cap B_r(x)$, $i=1,\ldots,n+1$, such that the $n$-dimensional ball
inscribed in the convex hull of $y_1, y_2, \ldots, y_{n+1}$, has radius not less than $c r$.
\end{theorem}
Smooth manifolds in $\R^n$ of dimension less than $n$ are examples of ``flat'' sets not preserving Markov's local inequality.

The interest to work with $d$-sets boundaries preserving  Markov's inequality (thus $0<d<n$), related in~\cite{BOS-1995} with Sobolev-Gagliardo-Nirenberg inequality, is to ensure the regular extensions $W^k_p(\Omega)\to W^k_p(\R^n)$ with $k\ge 2$ (actually the condition applies the continuity of the extension $C^\infty(\Omega)\to C^\infty(\R^n)$). For the extensions of minimal regularity $k=1$ 
(see in addition the Definition of Besov space Def.~3.2 in~\cite{ARXIV-IHNATSYEVA-2011}  with the help of the normalized local best approximation in the class of polynomials $P_{k-1}$ of the degree equal to $k-1$) Markov's inequality is trivially satisfied. 


Recently, Haj\l{}asz, Koskela and Tuominen~\cite{HAJLASZ-2008} have proved that every $W_p^k$-extension domain in $\R^n$ for $1\le p <\infty$ and $k\ge 1$, $k\in \N$ is an $n$-set. In addition, they proved that any $n$-set, for which $W_p^k(\Omega)=C_p^k(\Omega)$ (with norms' equivalence), is a $W_p^k$-extension domain for $1<p<\infty$ (see~\cite{HAJLASZ-2008} also for the results for $p=1$ and $p=\infty$). By $C_p^k(\Omega)$ is denoted the space of the fractional sharp maximal functions:
\begin{definition}
For a set $\Omega\subset \R^n$ of positive Lebesgue measure, $$C_p^k(\Omega)=\{f\in L_p(\Omega)|f_{k,\Omega}^\sharp(x)=\sup_{r>0} r^{-k}\inf_{P\in \mathcal{P}^{k-1}}\frac{1}{\mu(B_r(x))}\int_{B_r(x)\cap \Omega}|f-P|\dy\in L^p(\Omega)\}$$ with the norm $\|f\|_{C_p^k(\Omega)}=\|f\|_{L_p(\Omega)}+\|f_{k,\Omega}^\sharp\|_{L_p(\Omega)}.$
\end{definition}



From~\cite{JONES-1981} and~\cite{HAJLASZ-2008} we directly have
\begin{corollary}
 Let $\Omega$ be a bounded finitely connected domain in $\R^2$ and $1<p<\infty$, $k\in \N^*$. The domain $\Omega$ is a $2$-set with $W_p^k(\Omega)=C_p^k(\Omega)$ (with norms' equivalence) if and only if $\Omega$ is an $(\eps,\delta)$-domain and its boundary $\del \Omega$ consists of a finite number of points and quasi-circles. 
\end{corollary}

The question about $W^k_p$-extension domains is equivalent to the question of the continuity of the trace operator  $\mathrm{Tr}: W^k_p(\R^n) \to W^k_p(\Omega)$. Thus, let us  generalize the notion of the trace:
\begin{definition}\label{DefGTrace}
 For an arbitrary open set $\Omega$ of $\R^n$, the trace operator $\mathrm{Tr}$ is defined~\cite{JONSSON-1984,BODIN-2005,LANCIA-2002} for $u\in L_1^{loc}(\Omega)$ by
$$
 \mathrm{Tr} u(x)=\lim_{r\to 0} \frac{1}{\mu(\Omega\cap B_r(x))}\int_{\Omega\cap B_r(x)}u(y)dy.
$$
The trace operator $\mathrm{Tr}$ is considered for all $x\in\overline{\Omega}$ for which the limit exists.
\end{definition}

Using this trace definition it holds the trace theorem on closed $d$-sets~\cite{JONSSON-1984} Ch.VII and~\cite{WALLIN-1991} Proposition~4:
\begin{theorem}~\label{ThJWF}
   Let $F$ be a closed $d$-set preserving Markov’s
 local inequality.
 Then  if $0 < d < n$, $1 < p <\infty$, and $\beta = k - \frac{(n - d)}{p} > 0$, then the trace operator $\mathrm{Tr}: W^k_p(\R^n)\to B^{p,p}_\beta(F)$ is
bounded linear surjection
with a bounded right inverse $E: B^{p,p}_\beta(F)\to W^k_p(\R^n)$.
 \end{theorem}
 The definition of the Besov space  $B^{p,p}_{\beta}(F)$ on a closed
$d$-set $F$ can be found, for instance, in Ref.~\cite{JONSSON-1984} p.135 and
Ref.~\cite{WALLIN-1991}.
 
Hence, we introduce the notion of admissible domains:
\begin{definition}\textbf{(Admissible domain)}\label{DefAdmis}
 A domain $\Omega\subset \R^n$  is called admissible if it is an $n$-set, such that for $1<p<\infty$ and $k\in \N^*$ $W_p^k(\Omega)=C_p^k(\Omega)$  as sets with equivalent norms (hence, $\Omega$ is a $W_p^k$-extension domain), with a closed $d$-set boundary $\del \Omega$, $0<d<n$, 
 preserving local Markov's inequality.
\end{definition}

 Therefore, we summarize useful in the what follows results (see~\cite{JONSSON-1984,WALLIN-1991,JONSSON-1995,HAJLASZ-2008}) for the trace and the extension operators (see~\cite{SHVARTSMAN-2010} for more general results for the case $p>n$):
\begin{theorem}\label{ThContTrace}
 Let $1<p<\infty$, $k\in \N^*$ be fixed. Let $\Omega$ be an admissible domain in $\R^n$. 
Then, for $\beta=k-(n-d)/p>0$, the following trace operators (see Definition~\ref{DefGTrace})
\begin{enumerate}
 \item $\mathrm{Tr}: W_p^k(\R^n)\to B^{p,p}_{\beta}(\del \Omega)\subset L_p(\del \Omega)$,
 \item $\mathrm{Tr}_\Omega:W_p^k(\R^n)\to W_p^k(\Omega)$,
 \item $\mathrm{Tr}_{\del \Omega}:W_p^k(\Omega)\to B^{p,p}_{\beta}(\del \Omega)$
\end{enumerate}
are linear continuous and surjective with linear bounded  right inverse, $i.e.$ extension, operators $E: B^{p,p}_{\beta}(\del \Omega)\to W_p^k(\R^n)$, $E_\Omega: W_p^k(\Omega)\to W_p^k(\R^n)$ and $E_{\del \Omega}: B^{p,p}_{\beta}(\del \Omega)\to W_p^k(\Omega)$.
\end{theorem}
\textbf{Proof.}
 It is a corollary of results given in Refs.~\cite{JONSSON-1984,WALLIN-1991,JONSSON-1995,HAJLASZ-2008}. Indeed, if $\Omega$ is admissible, then by Theorem~\ref{ThJWF}, the trace operator $\mathrm{Tr}: W_p^k(\R^n)\to B^{p,p}_{\beta}(\del \Omega)\subset L_p(\del \Omega)$ is linear continuous and surjective with linear bounded  right inverse $E: B^{p,p}_{\beta}(\del \Omega)\to W_p^k(\R^n)$ (point 1). On the other hand, by~\cite{HAJLASZ-2008}, $\Omega$ is a $W_p^k$-extension domain and $\mathrm{Tr}_\Omega:W_p^k(\R^n)\to W_p^k(\Omega)$ and $E_\Omega: W_p^k(\Omega)\to W_p^k(\R^n)$ are linear continuous (point 2). Hence, the embeddings
 $$B^{p,p}_{\beta}(\del \Omega) \to W_p^k(\R^n)\to W_p^k(\Omega) \quad \hbox{and} \quad W_p^k(\Omega)\to W_p^k(\R^n)\to B^{p,p}_{\beta}(\del \Omega)$$
 are linear continuous (point 3).
$\Box$

Note that for $d=n-1$, one has $\beta=\frac{1}{2}$ and $
B_\frac{1}{2}^{2,2}(\del \Omega)=H^\frac{1}{2}(\del \Omega)$ as usual in the case of the classical results~\cite{LIONS-1972,MARSCHALL-1987} for Lipschitz boundaries $\del \Omega$.
In addition, for $u,\;v\in H^1(\Omega)$ with $\Delta u\in L_2(\Omega)$, the Green formula still holds in the framework of dual Besov spaces on a closed $d$-set boundary of $\Omega$ (see~\cite{LANCIA-2002} Theorem 4.15 for the von Koch case in $\R^2$):
\begin{proposition}\label{PropGreen}\textbf{(Green formula)}
Let $\Omega$ be an admissible domain in $\R^n$ ($n\ge 2$) with a $d$-set boundary $\del \Omega$ such that $n-2<d<n$. Then for all $u,\;v\in H^1(\Omega)$ with $\Delta u\in L_2(\Omega)$ it holds the Green formula
\begin{equation}\label{FracGreen}
  \int_\Omega v\Delta u\dx=\langle \frac{\del u}{\del \nu}, 
 \mathrm{Tr}v\rangle _{((B^{2,2}_{\beta}(\del \Omega))', B^{2,2}_{\beta}(\del \Omega))}-\int_\Omega \nabla v \nabla u \dx,
 \end{equation}
 where $\beta=1-(n-d)/2>0$ and the dual Besov space $(B^{2,2}_{\beta}(\del \Omega))'=B^{2,2}_{-\beta}(\del \Omega)$ is introduced in~\cite{JONSSON-1995}. 
 
 Equivalently, for an admissible domain $\Omega$ the normal derivative of $u\in H^1(\Omega)$ with $\Delta u\in L_2(\Omega)$ on the $d$-set boundary $\del \Omega$ with $n-2<d<n$ is defined by~Eq.~(\ref{FracGreen}) as a linear and continuous functional on $B^{2,2}_{\beta}(\del \Omega)$.
\end{proposition}
\textbf{Proof.}
The statement follows, thanks to Theorem~\ref{ThContTrace}, from the surjective property of the continuous trace operator $\mathrm{Tr}_{\del \Omega}:H^1(\Omega)\to B^{2,2}_{\beta}(\del \Omega)$ (see ~\cite{LANCIA-2002} Theorem 4.15). 
For a Lipschitz domain ($d=n-1$ and thus $
B_\frac{1}{2}^{2,2}(\del \Omega)=H^\frac{1}{2}(\del \Omega)$), we find the usual Green formula~\cite{LIONS-1972,MARSCHALL-1987}
\begin{equation*}
  \int_\Omega v\Delta u\dx=\langle \frac{\del u}{\del \nu}, 
 \mathrm{Tr}v\rangle _{((H^{\frac{1}{2}}(\del \Omega))', H^{\frac{1}{2}}(\del \Omega))}-\int_\Omega \nabla v \nabla u \dx.
 \end{equation*}
$\Box$

We also state the compact embedding of $H^1$ in $L_2$ for admissible truncated domains:
\begin{definition}\textbf{(Admissible truncated domain)}
 A domain $\Omega_S\subset \R^n$ ($n\ge 2$) is called admissible truncated domain  of an  exterior and admissible,  according to Definition~\ref{DefAdmis}, domain $\Omega$ with a compact $d_\Gamma$-set boundary $\Gamma$, if it is truncated by an admissible bounded domain $\Omega_1$ with a $d_S$-set boundary $S$,  $\Gamma\cap S=\varnothing$ (see Fig.~\ref{FigSchema}).
\end{definition}

 \begin{proposition}\label{PropCompEmb}
  Let 
  $\Omega_S$ be an admissible truncated domain with $n-2< d_S<n$.
  Then  the Sobolev space $H^1(\Omega_S)$ is compactly embedded in $L_2(\Omega_S)$: $$H^1(\Omega_S) \subset \subset L_2(\Omega_S).$$
 \end{proposition}
 \textbf{Proof.}
 Actually, in the case of a truncated domain, it is natural to impose $n-1\le d_\Gamma<n$ and $n-1\le d_S<n$, but formally the condition $\beta=1-(n-d)/2>0$ (see Theorem~\ref{ThContTrace} and Proposition~\ref{PropGreen}) only imposes the restriction  $n-2< d$.
 
 If $\Omega$ is an admissible domain (exterior or not), by Theorem~\ref{ThContTrace}, there exists linear bounded operator $E_{\Omega}: H^1(\Omega)\to H^1(\R^n)$.  
  Now, let in addition $\Omega$ be an exterior domain. Let us prove that for the admissible truncated domain $\Omega_S$ the extension operator $E_{\Omega_S\to \Omega}:H^1(\Omega_S) \to H^1(\Omega)$ is  a linear bounded operator.
 
 It follows from the fact that it is possible to extend $\Omega_1$ to $\R^n$ (there exists a linear bounded operator $E_{\Omega_1}: H^1(\Omega_1)\to H^1(\R^n)$) and that the properties of the extension are local, $i.e.$ depend on the properties of the boundary $S=\del \Omega_1$, which has no intersection with $\Gamma=\del \Omega$. For instance, if $S\in C^1$ (and thus $d_S=n-1$), then we can use the standard "reflection method" (as for instance in \cite{ALLAIRE-2012} Proposition 4.4.2). More precisely, we have to use a finite open covering $(\omega_i)_i$ of $S$ such that for all $i$ $\omega_i \cap \overline{\Omega_0} = \emptyset$. The compactness of $S$ and the fact that $S\cap \Gamma=\varnothing$ ensure that such a covering exists. In the case of a $d$-set boundary we use the Whitney extension method and Theorem~\ref{ThJWF}. 
  
  Hence, using Theorem~\ref{ThContTrace}, there exists a linear bounded operator $A: H^1(\Omega)\to H^1(\Omega_1)$ as a composition of the extension operator $E_{\Omega}: H^1(\Omega)\to H^1(\R^n)$ and the trace operator $Tr_{\Omega_1}:H^1(\R^n)\to H^1(\Omega_1)$ ($A=Tr_{\Omega_1}\circ E_{\Omega}$).
  
  Let us define a parallelepiped $\Pi$ in the such way that
$$\Omega_1 \subsetneq \Pi,\quad \Pi=\{x=(x_1,\ldots,x_n)|\;0<x_i<d_i \; (d_i\in \R)\}.$$
Consequently, the operator
  $B=E_{\Omega_1\to \Pi}\circ  A\circ E_{\Omega_S\to \Omega} :H^1(\Omega_S) \to H^1(\Pi)$ is a linear bounded operator as a composition of linear bounded operators.  

  Let us prove $H^1(\Omega_S)\subset \subset L_2(\Omega_S)$. We follow  the proof of the compact embedding of $H^1$ to $L_2$, given in~\cite{MASLENNIKOVA-1997} in the case of a regular boundary.
  
  Indeed, let $(u_m)_{m\in\N}$ be a bounded sequence in $H^1(\Omega_S)$. 
Thanks to the boundness of the operator $B$, for all $m\in \N$ we  extend $u_m$ from $\Omega_S$ to the parallelepiped $\Pi$, containing $\Omega_S$. Thus, for all $m\in \N$ the extensions $Bu_m=\tilde{u}_m$ satisfy
$$\tilde{u}_m\in H^1(\Pi),\quad \tilde{u}_m|_{\Omega_S}=u_m,\quad \|u_m\|_{H^1(\Omega_S)}\le \|\tilde{u}_m\|_{H^1(\Pi)},$$
and, in addition, there exists a constant $C(\Omega,\Pi)=\|B\|>0,$ independent on $u_m$, such that
$$\|\tilde{u}_m\|_{H^1(\Pi)}\le C(\Omega,\Pi) \|u_m\|_{H^1(\Omega_S)}.$$
Hence, the sequence $(\tilde{u}_m)_{m\in\N}$ is also a  bounded sequence in $H^1(\Pi)$. Since the embedding $H^1(\Pi)$ to $L_2(\Pi)$ is continuous, the sequence $(\tilde{u}_m)_{m\in\N}$ is also bounded in $L_2(\Pi)$.
Let in addition $$\Pi=\sqcup_{i=1}^{N^n} \Pi_i, \hbox{ where } \Pi_i=\otimes_{k=1}^n [a_i,a_i+\frac{d_k}{N}] \quad (a_i\in \R).$$
Thanks to~\cite{MASLENNIKOVA-1997} p.~283, in  $\Pi$ there holds the following inequality for all $u\in H^1(\Pi)$:
\begin{equation}\label{eqeq}
 \int_\Pi u^2 dx\le \sum_{i=1}^{N^n}\frac{1}{|\Pi_i|}\left(\int_{\Pi_i}udx\right)^2+\frac{n}{2}\int_\Pi \sum_{k=1}^n\left(\frac{d_k}{N}\right)^2\left(\frac{\del u}{\del x_k} \right)^2dx.
\end{equation}

On the other hand,  $L_2(\Pi)$ is a  Hilbert space, thus weak$^*$ topology on it is equal to the weak topology. Moreover, as $L_2$ is separable,  all closed bounded sets in $L_2(\Pi)$ are weakly sequentially compact (or compact in the weak topology since here the weak topology is metrizable). 
To simplify the notations,  we  simply write $u_m$ for $\tilde{u}_m \in L_2(\Pi)$.
Consequently, the sequence $(u_m)_{m\in\N}$ is weakly sequentially compact in $L_2(\Pi)$ and we have
$$\exists (u_{m_k})_{k\in\N}\subset (u_m)_{m\in\N} \;: \quad \exists u\in L_2(\Pi) \quad u_{m_k}\rightharpoonup u \hbox{ for } k\to+\infty.$$
Here $u$ is an element of $L_2(\Pi)$, not necessarily in $H^1(\Pi)$.

As $(L_2(\Pi))^*=L_2(\Pi)$, by the Riesz representation theorem, 
$$u_{m_k}\rightharpoonup u \in L_2(\Pi) \quad \Leftrightarrow \quad 
\forall v \in L_2(\Pi) \quad \int_\Pi(u_{m_k}-u)vdx\to 0.$$
Since $(u_{m_k})_{k\in\N}$ is a Cauchy sequence in the weak topology on $L_2(\Pi)$, then, in particular choosing $v=\mathds{1}_\Pi$, it holds $$\int_\Pi(u_{m_k}-u_{m_j})dx\to 0 \quad \hbox{for } k,j\to +\infty.$$
 
Thus, using Eq.~(\ref{eqeq}), for two members of the sub-sequence $(u_{m_k})_{k\in\N}$ with sufficiently large ranks $p$ and $q$, we have
\begin{multline*}
 \|u_p-u_q\|^2_{L_2(\Omega_S)}\le \|u_p-u_q\|^2_{L_2(\Pi)}\\
 \le \sum_{i=1}^{N^n}\frac{1}{|\Pi_i|}\left(\int_{\Pi_i}(u_p-u_q)dx\right)^2+\frac{n}{2N^2}\sum_{k=1}^nd_k^2\left\| \frac{\del u_p}{\del x_k}-\frac{\del u_q}{\del x_k}\right\|^2_{L_2(\Pi)}<\frac{\eps}{2}+\frac{\eps}{2}=\eps.
\end{multline*}
Here we have chosen $N$ such that $$\frac{n}{2N^2}\sum_{k=1}^nd_k^2\left\| \frac{\del u_p}{\del x_k}-\frac{\del u_q}{\del x_k}\right\|^2_{L_2(\Pi)}<\frac{\eps}{2}.$$
Consequently, $(u_{m_k})_{k\in\N}$ is a Cauchy sequence in $L_2(\Omega_S)$, and thus converges strongly in $L_2(\Omega_S)$. 
 $\Box$
 \begin{remark}
To have a compact embedding  it is important  that 
 the domain $\Omega$ be an $W^k_p$-extension domain.
The boundness or unboudness of $\Omega$ is not important to have $W_p^{k}(\Omega)\subset \subset W_p^{\ell}(\Omega)$ with $k>\ell\ge1$ ($1<p<\infty$).
But the boundness of $\Omega$ is important for the  compact embedding in $L_q(\Omega)$.
\end{remark}
As a direct corollary we have the following generalization of the classical Rellich-Kondrachov theorem (see for instance Adams~\cite{ADAMS-1975} p.144 Theorem 6.2):
\begin{theorem}\label{ThCSEnSET} \textbf{(Compact Sobolev embeddings for $n$-sets)}
 Let $\Omega\subset \R^n$ be 
 an $n$-set with $W_p^k(\Omega)=C_p^k(\Omega)$, $1<p<\infty$, $k,\ell\in \N^*$. Then there hold  the following compact  embeddings:    
    \begin{enumerate}
     \item $W_p^{k+\ell}(\Omega)\subset \subset W_q^\ell(\Omega)$,
     \item $W_p^{k}(\Omega)\subset \subset L_q^{loc}(\Omega)$,  or $W_p^{k}(\Omega)\subset \subset L_q(\Omega)$ if $\Omega$ is bounded,
    \end{enumerate}
 with $q\in[1,+\infty[$ if $kp=n$, $q\in [1,+\infty]$ if $kp>n$, and with $q\in[1, \frac{pn}{n-kp}[$ if $kp<n$. 
\end{theorem}
\textbf{Proof.}
Let us denote by  $B_r(x)$ a non trivial ball for the Euclidean metric in $\R^n$ (its boundary is infinitely smooth, and thus, it is a $W^k_p$-extension domain for all  $1<p<\infty$ and  $k\in \N^*$).  By~\cite{HAJLASZ-2008} (see also Theorem~\ref{ThContTrace}),  the extension $E: W_p^{k+\ell}(\Omega)\to W_p^{k+\ell}(\R^n)$ and the trace $\mathrm{Tr}_{B_r}:W_p^{k+\ell}(\R^n)\to W_p^{k+\ell}(B_r(x))$ are continuous. In addition, by  the classical Rellich-Kondrachov theorem on the ball $B_r(x)$, the embedding $K: W_p^{k+\ell}(B_r(x))\to W_q^\ell(B_r(x))$, for the mentioned values of $k$, $p$, $n$ and $\ell$, is compact.
 Hence, for $\ell\ge 1$,  thanks again to~\cite{HAJLASZ-2008}, $E_1: W_q^\ell(B_r(x)) \to W_q^\ell(\Omega)$ is continuous, as the composition of continuous operators  $E_2: W_q^\ell(B_r(x)) \to W_q^\ell(\R^n)$ and $\mathrm{Tr}_{\Omega}:W_q^\ell(\R^n) \to W_q^\ell(\Omega)$. Finally, the embedding $W_p^{k+\ell}(\Omega)\subset W_q^\ell(\Omega)$ for $\ell\ge 1$ is compact, by the composition of the continuous and  compact operators:
 $$E_1\circ K \circ \mathrm{Tr}_{B_r} \circ E: W_p^{k+\ell}(\Omega)\to W_q^\ell(\Omega).$$
 When $\ell=0$, instead of Sobolev embedding $E_1$, we need to have the continuous embedding of Lebesgue spaces
 $L_q(B_r(x)) \to L_q(\Omega)$, which holds if and only if $\Omega$ is bounded. If $\Omega$ is not bounded, for all measurable compact sets $K\subset \Omega$, the embedding $L_q(B_r(x))\to L_q(K)$ is continuous. This finishes the proof.
 $\Box$

In the same way, we generalize  the classical Rellich-Kondrachov theorem for fractals: 
\begin{theorem}\label{ThAnnaEmbF} \textbf{(Compact Besov embeddings for $d$-sets)}
 Let $F\subset \R^n$ be a closed $d$-set preserving Markov's local inequality, $0<d<n$, $1<p<\infty$ and $\beta=k+\ell-\frac{n-d}{p}>0$ for $k,\ell\in \N^*$. 
 
 Then, for the same $q$ as in Theorem~\ref{ThCSEnSET}, 
 the following continuous embeddings  are compact
 \begin{enumerate}
  \item $B^{p,p}_{\beta}(F)\subset \subset B^{q,q}_\alpha(F)$ for  $\ell\ge 1$ and $\alpha=\ell-\frac{n-d}{q}>0$;
  \item if $F$ is bounded in $\R^n$, $B^{p,p}_{\beta}(F)\subset \subset L_q(F)$, otherwise $B^{p,p}_{\beta}(F)\subset \subset L_q^{loc}(F)$ for $\ell\ge 0$.
 \end{enumerate}
 \end{theorem}
\textbf{Proof.}
Indeed, thanks to Theorem~\ref{ThJWF}, the extension
$E_F:B^{p,p}_{\beta}(F) \to W_p^{k+\ell}(\R^n)$ is continuous. Hence, by Calderon~\cite{CALDERON-1961}, a non trivial ball is $W_p^{k+\ell}$-extension domain:
$\mathrm{Tr}_{B_r}$ (see the proof of Theorem~\ref{ThCSEnSET}) is continuous. Thus, the classical Rellich-Kondrachov theorem on the ball $B_r(x)$ gives the compactness of $K: W_p^{k+\ell}(B_R)\to W_q^\ell(B_r)$. Since, for $\ell\ge 1$, $E_2:W_q^\ell(B_r)\to W_q^\ell(\R^n)$ is continuous and, by Theorem~\ref{ThJWF},
 $\mathrm{Tr}_F:W_q^\ell(\R^n)\to B^{q,q}_\alpha(F)$ is continuous too, we conclude that the operator $$\mathrm{Tr}_F\circ K \circ \mathrm{Tr}_{B_r} \circ E_F:B^{p,p}_{\beta}(F) \to B^{q,q}_\alpha(F)$$
 is compact. For $j=0$, we have $W_q^0=L_q$, and hence, if $F\subset B_r(x)$, the operator 
 $ L_q(B_r(x)) \to L_q(F)$ is a linear continuous measure-restriction operator on a $d$-set (see~\cite{JONSSON-1984} for the $d$-measures). If $F$ is not bounded in $\R^n$,  for all bounded $d$-measurable subsets $K$ of $F$, the embedding $L_q(B_r(x))\to L_q(K)$ is continuous. 
 $\Box$
 
In particular, the compactness of the trace operator implies the following equivalence of the norms on $W_p^k(\Omega)$:
 \begin{proposition}\label{PropCompET}
 Let $\Omega$ be an admissible domain in $\R^n$ with a compact boundary $\del \Omega$ and $1<p<\infty$, $k\in \N^*$, $\beta=k-\frac{n-d}{p}>0$. Then
 \begin{enumerate}
  \item $W_p^k(\Omega)\subset \subset L_p^{loc}(\Omega)$;
  \item $\mathrm{Tr}: W_p^k(\Omega)\to L_p(\del \Omega)$ is compact;
  \item $\|u\|_{W_p^k(\Omega)}$ is equivalent to $\|u\|_{\mathrm{Tr}}=\left(\sum_{|l|=1}^k \int_\Omega |D^l u|^p\dx +\int_{\del \Omega} |\mathrm{Tr}u|^p\dm_d \right)^\frac{1}{p}. $
 \end{enumerate}
 
\end{proposition}
\textbf{Proof.}
Point~1 follows from Theorem~\ref{ThCSEnSET} and holds independently on values of $kp$ and $n$.
The trace operator $\mathrm{Tr}: W_p^k(\Omega)\to L_p(\del \Omega)$ from  Point~2 is compact as a composition of the compact, by Theorem~\ref{ThAnnaEmbF}, operator $K:B^{p,p}_{\beta}(\del \Omega)\to L_p(\del \Omega)$ with the continuous operator $\mathrm{Tr}_{\del \Omega}: W_p^k(\Omega)\to B^{p,p}_{\beta}(\del \Omega)$.

Let us show that Points 1 and 2 imply the equivalence of the norms in Point~3.
 We generalize the proof of Lemma 2.2 in Ref.~\cite{ARENDT-2015}.
 Since the trace $\mathrm{Tr}: W_p^k(\Omega)\to L_p(\del \Omega)$ is continuous, then there exists a constant $C>0$ such that for all $u\in W_p^k(\Omega)$
 $$\int_{\del \Omega} |\mathrm{Tr}u|^p\dm_d\le \sum_{|l|=1}^k \int_\Omega |D^l u|^p\dx +\int_{\del \Omega} |\mathrm{Tr}u|^p\dm_d\le C\|u\|_{W_p^k(\Omega)}^p.$$
 Let us prove that there a constant $c>0$ such that for all $u\in W_p^k(\Omega)$
 $$\int_{ \Omega} |u|^p\dx\le c\left(\sum_{|l|=1}^k \int_\Omega |D^l u|^p\dx +\int_{\del \Omega} |\mathrm{Tr}u|^p\dm_d\right)= c\|u\|_{\mathrm{Tr}}^p.$$
 Suppose the converse. Then for all $m\in \N^*$ there exists a $u_m\in W_p^k(\Omega)$ such that
 \begin{equation}\label{EqConverse}
  \|u_m\|^p_{\mathrm{Tr}}<\frac{1}{m}\int_{ \Omega} |u_m|^p\dx.
 \end{equation}
 As in Ref.~\cite{ARENDT-2015}, without loss of generality we assume that $$\hbox{for all }m\in \N^*\quad \|u_m\|_{L_p(\Omega)}=1.$$
 Then the sequence $(u_m)_{m\in\N}$ is bounded in $W_p^k(\Omega)$: for all $m\in \N^*$ $\|u_m\|^p_{W_p^k(\Omega)}\le 2$. As $W_p^k(\Omega)$ is a reflexive Banach space, each bounded sequence in  $W_p^k(\Omega)$ contains a weakly convergent subsequence. Hence, there exists $u\in  W_p^k(\Omega)$ such that $u_{m_i}\rightharpoonup u$ in $W_p^k(\Omega)$ for $m_i \to +\infty$. By the compact embedding of $W_p^k(\Omega)$ in $L_p(\Omega)$ (Point 1), the subsequence $(u_{m_i})_{i\in\N}$ converges strongly towards $u$ in $L_p(\Omega)$. Consequently, $\|u\|_{L_p(\Omega)}=1$ and
 $$\sum_{|l|=1}^k \int_\Omega |D^l u|^p\dx\le \liminf_{i\to +\infty}  \sum_{|l|=1}^k \int_\Omega |D^l u_{m_i}|^p\dx\le \liminf_{i\to +\infty} \frac{1}{n}=0.$$
 Therefore, $u$ is constant (since $\Omega$ is connected) with $\|u\|_{L_p(\Omega)}=1$. From Eq.~(\ref{EqConverse}) we have
 $$
 \int_{\del \Omega} |\mathrm{Tr}u_{m_i}|^p\dm_d\le \|u_{m_i}\|^p_{\mathrm{Tr}}<\frac{1}{m_i}\int_{ \Omega} |u_{m_i}|^p\dx.$$
 
 Since the trace operator $\mathrm{Tr}: W_p^k(\Omega)\to L_p(\del \Omega)$ is compact (Point 2), it holds
 $$\|\mathrm{Tr}u\|^p_{L_p(\del \Omega)}=\lim_{i\to +\infty} \|\mathrm{Tr}u_{m_i}\|^p_{L_p(\del \Omega)}=0,$$
 which implies that $u=0$. This is a contradiction with $\|u\|_{L_p(\Omega)}=1$. Hence, there exists a constant $\tilde{c}>0$ such that
 $\|u\|_{W_p^k(\Omega)}\le \tilde{c}\|u\|_{\mathrm{Tr}}$.
$\Box$

Going back to the Laplace transport problem on  exterior and truncated domains, we especially need the following theorem, thus formulated for $H^1$:
\begin{theorem}\label{ThCompTrExtdom}\textbf{(Compactness of the trace)}
 Let $\Omega$ be an admissible domain of $\R^n$ with a compact $d$-set boundary $\Gamma$, $n-2< d_\Gamma<n$. If $\Omega$ is an exterior domain, let $\Omega_S$ be its admissible truncated domain with $n-2< d_S<n$.
  Then for all these domains, $i.e.$ for $D=\Omega$ or $D=\Omega_S$,
 there exist linear trace operators 
 $$\mathrm{Tr}_\Gamma:  H^1(D) \rightarrow L_2(\Gamma)\quad\hbox{ and, if } S\ne \varnothing, \quad \mathrm{Tr}_S:  H^1(\Omega_S) \rightarrow L_2(S),$$
 which are compact. Moreover, $\mathrm{Im(Tr_\Gamma}(H^1(D)))=B^{2,2}_{\beta_\Gamma}(\Gamma)$ for $\beta_\Gamma=1-\frac{n-d_\Gamma}{2}>0$ and 
 $\mathrm{Im(Tr_S}(H^1(\Omega_S)))=B^{2,2}_{\beta_S}(S)$ for $\beta_S=1-\frac{n-d_S}{2}>0$.
 \end{theorem}
\textbf{Proof.}
It is a direct corollary of Proposition~\ref{PropCompET}.
$\Box$


 \section{Well-posedness  of Robin boundary problem for the Laplace equation}\label{SecWellPosed}
        \subsection{Well-posedness on a truncated domain}\label{SecWellPosedTrunc}
  Let us start by a well-posedness of  problem~(\ref{Problem}) for an admissible truncated domain $\Omega_S$ introduced in Section~\ref{SecIntro}.
         Therefore,  
      $\Omega_S$ is a bounded domain with a compact $d_\Gamma$-set boundary $\Gamma$, $n-2< d_\Gamma<n$ ($n\ge 2$), on which is imposed  the Robin boundary condition for $\lambda>0$ and $\psi\in L_2(\Gamma)$, and a  $d_S$-set boundary $S$, $n-2< d_S<n$, on which is imposed  the homogeneous Dirichlet boundary condition.

    Let us denote  $\tilde{H}^1(\Omega_S) := \{ u \in H^1(\Omega_S) : \mathrm{Tr}_S u = 0 \}$. Note that, thanks to Theorem~\ref{ThCompTrExtdom},
    the continuity of the map $Tr_S$ ensures  that $\tilde{H}^1(\Omega_S)$ is a Hilbert space.
Therefore, thanks to Proposition~\ref{PropCompEmb}, as $H^1(\Omega_S)\subset \subset L_2(\Omega_S)$, following for instance the proof of Evans~\cite{EVANS-1994} (see section 5.8.1 Theorem 1), we obtain
\begin{proposition}
 \label{PropPoincare} \textbf{(Poincar\'e inequality)}
 Let $\Omega_S\subset \R^n$ be an admissible truncated domain, introduced in Theorem~\ref{ThCompTrExtdom}, with $n\ge 2$.
 For all  $v \in \tilde{H}^1(\Omega_S)$ there exists $C > 0$, depending only on $\Omega_S$ and $n$, such that
\[
\| v \|_{L_2(\Omega_S)} \leq C \| \nabla v \|_{L_2(\Omega_S)}.
\]
 Therefore the semi-norm $\|\cdot\|_{\tilde{H}^1(\Omega_S)}$, defined by $\|v\|_{\tilde{H}^1(\Omega_S)}=\|\nabla v \|_{L_2(\Omega_S)}$, is a norm, which is equivalent to $\|\cdot\|_{H^1(\Omega_S)}$ on $\tilde{H}^1(\Omega_S)$.
\end{proposition}
\begin{remark}\label{RemPoincGen}
 Let us denote $\langle v \rangle=\frac{1}{Vol(\Omega_S)} \int_{\Omega_S} v \dx$. Since $\Omega_S$ is a bounded $W^1_p$-extension domain,   Theorem~\ref{ThCSEnSET}  ensures $W^1_p(\Omega_S)\subset\subset L_p(\Omega_S)$ for all $1< p< \infty$. Thus the  Poincar\'e inequality can be generalized with the same proof to $W^1_p(\Omega_S)$ for all $1< p< \infty$:
 $$\forall v \in W^1_p(\Omega_S) \; \exists C=C(\Omega_S,p,n)>0:\quad \| v-\langle v \rangle \|_{L^p(\Omega_S)} \leq C \| \nabla v \|_{L^p(\Omega_S)}.$$
\end{remark}

Consequently, using these properties of $\tilde{H}^1(\Omega_S)$, we have the well-posedness of problem~(\ref{Problem}):

     \begin{theorem}\label{ThWelPTr}
     Let $\Omega_S\subset \R^n$ be an admissible truncated domain, introduced in Theorem~\ref{ThCompTrExtdom}, with $n\ge 2$.
      For all $\psi \in L_2(\Gamma)$ and $\lambda \geq 0$, there exists a unique weak solution $u\in \tilde{H}^1(\Omega_S)$ of problem~(\ref{Problem}) such that
      \begin{equation}\label{EqVarFTruncD}
       \forall v \in \tilde{H}^1(\Omega_S) \quad \int_{\Omega_S} \nabla u \nabla v \dx + \lambda\int_\Gamma \mathrm{Tr}_\Gamma u \mathrm{Tr}_\Gamma v \mathrm{d} m_{d_\Gamma} = \int_\Gamma \psi \mathrm{Tr}_\Gamma v \mathrm{d} m_{d_\Gamma}.
      \end{equation}
Therefore, for all $\lambda\in [0,\infty[$ and $\psi\in L_2(\Gamma)$ the operator $$B_\lambda(S): \psi \in L_2(\Gamma)\mapsto u\in \tilde{H}^1(\Omega_S)$$ with $u$, the solution of the variational problem~(\ref{EqVarFTruncD}), has the following properties
    \begin{enumerate}
     \item $B_\lambda(S)$ is a linear compact operator;
     \item $B_\lambda(S)$ is positive: if $\psi\ge 0$ from $L_2(\Gamma)$, then  for all $\lambda \in [0, \infty[$ $B_\lambda \psi=u\ge 0$;
     \item $B_\lambda(S)$ is monotone: if $0\le \lambda_1<\lambda_2$, then for all $\psi\ge 0$ from $L_2(\Gamma)$ it holds $B_{\lambda_2}(S) \psi=u_{\lambda_2} \le u_{\lambda_1}= B_{\lambda_1}(S) \psi$;
     \item If $\lambda \in [0, \infty[$ then $0\le B_\lambda(S) \mathds{1}_\Gamma\le \frac{1}{\lambda} \mathds{1}_{\Omega_S}$.
         \end{enumerate}   
         \end{theorem}
 \textbf{Proof.}
  It's a straightforward application of the Lax-Milgram theorem. The continuity of the two forms is ensured by the continuity of the trace operator $\mathrm{Tr}_\Gamma$ (see   Theorem~\ref{ThCompTrExtdom}). The coercivity of the symmetric bilinear form is ensured by the Poincar\'e inequality (see Proposition~\ref{PropPoincare}). To prove the properties of the operator $B_\lambda(S)$ it is sufficient to replace $W^D(\Omega)$ by  $\tilde{H}^1(\Omega_S)$ in the proof of Theorem~\ref{ThweakFormWellPosed}. 
 $\Box$

        \subsection{Functional spaces for the exterior problem}~\label{SubsFuncSpW}
          To be able to prove the well-posedness of problem~(\ref{Problem}) on an exterior domain with Dirichlet boundary conditions at infinity, we  extend the notion of $(\tilde{H}^1,\|\cdot\|_{\tilde{H}^1})$ to the unbounded domains. If $\Omega$ is an exterior domain of a bounded domain $\Omega_0$, $i.e.$ $\Omega=\R^n\setminus \overline{\Omega}_0$, the usual Poincar\'e inequality does not hold anymore and, hence, we don't have Proposition~\ref{PropPoincare}.
          For this purpose, we use~\cite{LU-2005,ARENDT-2015} and define for $\Omega=\R^n\setminus \overline{\Omega}_0$, satisfying the conditions of Theorem~\ref{ThContTrace},
          \begin{equation*}
          W(\Omega) := \{ u \in H^1_{loc}(\Omega), \quad \int_\Omega |\nabla u|^2\dx < \infty \}.
          \end{equation*}
          \begin{remark}
            Let us fix a $r_0>0$ in the way  that there exists $x\in \R^n$ such that $\overline{\Omega}_0\subset B_{r_0}(x)=\{y\in \R^n| \;|x-y|<r_0\}$, and for all $r\ge r_0$ define $\Omega_r=\Omega\cap B_r(x)$. Thanks to Remark~\ref{RemPoincGen}, locally we always have the Poincar\'e inequality: 
          \begin{equation*}
           \forall u\in W(\Omega) \quad \|u-\langle u \rangle\|_{L_2^{loc}(\Omega)}\le C_{loc}\|\nabla u\|_{L_2^{loc}(\Omega)}\le C_{loc}\|\nabla u\|_{L_2(\Omega)}<\infty,
          \end{equation*}
which implies that for all $r\ge r_0$ the trace $u|_{\Omega_r}\in H^1(\Omega_r)$ (see Proposition~\ref{PropCompET} and Theorem~\ref{ThCompTrExtdom}).
Therefore, as in~\cite{ARENDT-2015}, it is still possible to consider (but we don't need it)
\begin{equation*}
          W(\Omega) = \{ u : \Omega \rightarrow \C | u \mbox{ is measurable, } \forall r > r_0\; u|_{\Omega_r} \in H^1(\Omega_r)  \mbox{ and }\int_\Omega |\nabla u|^2 < \infty \}.
          \end{equation*}
          \end{remark}          
          
Thanks to G. Lu and B. Ou (see \cite{LU-2005} \textbf{Theorem 1.1} with $p=2$), we have
\begin{theorem}
\label{theorem/Ou1}
Let $u \in W(\R^n)$ with $n\ge 3$. Then there exists the following limit:
\begin{eqnarray*}
(u)_{\infty} = \lim_{R \rightarrow \infty} \frac{1}{| B_R |} \int_{B_R} u\dx.
\end{eqnarray*}
Moreover, there exists a constant $c > 0$, depending only on the dimension $n$, but not on $u$, such that:
\begin{equation*}\label{IneqvLL}
 \| u - (u)_{\infty} \|_{L_{\frac{2n}{n - 2}}(\R^n)} \leq c \| \nabla u \|_{L_2(\R^n)}.
\end{equation*}
\end{theorem}

In~\cite{LU-2005} Section 5, G. Lu and B. Ou extend this result to exterior domains with a Lipschitz boundary. Their proof is based on the existence of a continuous extension operator. Therefore, thanks to Theorem~\ref{ThContTrace}, we generalize  Theorem 5.2 and Theorem 5.3 of G. Lu and B. Ou and take $p=2$, according to our case.
\begin{theorem}
\label{theorem/Ou2}
Let $n\ge 3$ and $\Omega$ be an admissible exterior domain  with a compact $d$-set boundary $\Gamma$ ($n-2< d<n$). 
There exists $c := c(n, \Omega) > 0$ so that for all $u \in W(\Omega)$ there exists $(u)_{\infty} \in \R$ such that 
\begin{equation}\label{InOu}
 (\int_\Omega | u - (u)_{\infty}|^{\frac{2n}{n - 2}})^{\frac{n - 2}{2n}} \leq c(n, \Omega) \, \| \nabla u \|_{L_2(\Omega)}.
\end{equation}
Moreover, it holds
\begin{enumerate}
  \item The space $W(\Omega)$ is a Hilbert space, corresponding the inner product
 $$(u, v) := \int_\Omega \nabla u . \nabla v \dx+ (u)_{\infty}(v)_{\infty}.$$
 The associated norm is denoted by $\|u\|_{W(\Omega)}$.
 \item The following norms are equivalent to $\|\cdot\|_{W(\Omega)}$: 
 $$\|u\|_{\Gamma,\Omega}=(\|\nabla u\|^2_{L_2(\Omega)} + \|\mathrm{Tr} u\|^2_{L_2(\Gamma)} )^{1/2},\quad \|u\|_{A,\Omega}=( \|\nabla u\|^2_{L_2(\Omega)} + \|u\|^2_{L_2(A)} )^{1/2},$$
 where $A \subset \Omega$ is a bounded measurable set with $\mathrm{Vol}(A)=\int_A 1\dx >0$.
 \item There exists a continuous extension operator $E : W(\Omega) \rightarrow W(\R^d)$.
 \item The trace operator $\mathrm{Tr} : W(\Omega) \rightarrow L_2(\Gamma)$ is compact.
\end{enumerate}
\end{theorem}
\textbf{Proof.}
 Thanks to Theorem~\ref{ThContTrace}, we update Theorem 5.2 and 5.3~\cite{LU-2005} to obtain  inequality~(\ref{InOu}).
 Let us notice the importance of the Sobolev embedding $H^1(\R^n)\subset L_{\frac{2n}{n - 2}}(\R^n)$ which holds for $n\ge 3$, but which is false for $n=2$. The real number $(u)_{\infty}$ in  the inequality~(\ref{InOu}) is merely the 'average' of an extension of $u$ to $\R^n$, as defined in Theorem~\ref{theorem/Ou1}.
 
 Point 1, stating the completeness of $W(\Omega)$, follows from Ref.~\cite{LU-2005} by updating the proof of Theorem 2.1. 
 
   The equivalence of norms in Point~2 follows from the proof of Proposition~\ref{PropCompET} using Theorems~\ref{ThContTrace} and~\ref{ThCompTrExtdom} (see also Proposition 2.5~\cite{ARENDT-2015}). 
 
 To prove Point~3,
 we notice that, thanks to Point 2, the extension operator $E$ is continuous if and only if the domain $\Omega$ is such that the extension $E_\Omega: H^1(\Omega)\to H^1(\R^n)$ is a linear continuous operator. This is true in our case, since the domain $\Omega$ satisfies the conditions of Theorem~\ref{ThCompTrExtdom}.  
 In addition, the continuity of $E_\Omega$ ensures that, independently on the geometric properties of the truncated boundary $S$ ($S\cap \Gamma=\varnothing$), for all (bounded) truncated domains $\Omega_S$ the extension operator $E_0: H^1(\Omega_S)\to H^1(\Omega_S\cup \overline{\Omega}_0)$ is continuous. Indeed, if $E_\Omega: H^1(\Omega)\to H^1(\R^n)$ is continuous, then $H^1_{loc}(\Omega)\to H^1(\R^n)$ is also continuous and hence, 
  we can consider only functions with a support on $\Omega_S$ and extend them to $\Omega_S\cup \overline{\Omega}_0=\R^n\cap \Omega_1$ to obtain the continuity of $E_0$.
 

To prove Point~4, we write $\mathrm{Tr} : W(\Omega) \rightarrow L_2(\Gamma)$ as a composition of two traces operators:
\begin{multline*}
 Tr= Tr_\Gamma \circ Tr_{W\to H^1}, \quad Tr_{W\to H^1}: W(\Omega)\to H^1(\Omega_S),\quad
 Tr_\Gamma : H^1(\Omega_S)\to L_2(\Gamma).
\end{multline*}
As $Tr_{W\to H^1}$ is continuous, $i.e.$ $$\|u\|_{H^1(\Omega_S)}^2\le C(\|\nabla u\|_{L_2(\Omega)}^2+\| u\|_{L_2(\Omega_S)}^2),$$
and, since $\Omega$ is an admissible domain with a compact boundary $\Gamma$, by Proposition~\ref{PropCompET},  $Tr_\Gamma$ is compact, we deduce the compactness of $\mathrm{Tr} : W(\Omega) \rightarrow L_2(\Gamma)$. 
 $\Box$

To have an analogy in the unbounded case with $\tilde{H}^1$ for a truncated domain, let us introduce, as in~\cite{ARENDT-2015}, the space
$W^D(\Omega)$, defined by the closure of the space $$\{\,u|_{\Omega} \,:\, u \in \mathcal{D}(\R^n),\, \, n\ge 3\}$$ with respect to the norm $u \mapsto (\int_\Omega |\nabla u|^2)^{1/2}$. Therefore, for  the inner product 
$$
(u, v)_{W^D(\Omega)} = \int_\Omega \nabla\,u . \nabla\,v
,$$ the space $(W^D(\Omega), (., .)_{W^D(\Omega)})$ is a Hilbert space (see a discussion about it on p.~8 of Ref.~\cite{LU-2005}).


It turns out that $W^D(\Omega)$ is the space of all $u \in W(\Omega)$ with average zero: 
\begin{proposition}
Let $\Omega$ be a unbounded (actually, exterior) domain in $\R^n$ with $n\ge 3$. The space $W^D(\Omega)$ has co-dimension 1 in $W(\Omega)$. Moreover $$W^D(\Omega) = W(\Omega) \cap L_{\frac{2n}{n - 2}}(\Omega) = \{u \in W(\Omega) : (u)_\infty = 0\}.$$
\end{proposition}
\textbf{Proof.}
 See~\cite{ARENDT-2015} Proposition 2.6 and references therein.
 $\Box$
 \begin{remark}
  Note that, as $n\ge 3$, $H^1(\Omega) \subset W(\Omega) \cap L_{\frac{2n}{n - 2}}(\Omega)= W^D(\Omega)$, which is false for $n=2$.
 \end{remark}

        \subsection{Well-posedness of the exterior problem and its approximation}
        

       Given $\psi \in L_2(\Gamma)$ and $\lambda \geq 0$, we consider the Dirichlet problem on the exterior domain $\Omega$ with Robin boundary conditions on the boundary $\Gamma$ in $\R^n$, $n\ge 3$: 
\begin{eqnarray}
\label{eq/RobCond}
-\Delta u &=& 0 \quad x \in \Omega, \nonumber\\
\lambda \mathrm{Tr} u + \partial_\nu u &=& \psi \quad x \in \Gamma. 
\end{eqnarray}
      At infinity we consider Dirichlet boundary conditions. In \cite{ARENDT-2015} W. Arendt and A.F.M ter Elst also considered Neumann boundary conditions at infinity. Those results apply as well in our setting, but we chose to focus on the Dirichlet boundary conditions at infinity in order not to clutter the presentation. It is worth emphasizing that in the following we only consider weak formulations that we describe below. 

      Since (see Subsection~\ref{SubsFuncSpW}) $H^1(\Omega)\subset W^D(\Omega)\subset W(\Omega)$ by their definitions,  we need to update the definition of the normal derivative, given by Eq.~(\ref{FracGreen}) in Section~\ref{SecFuncAnT}, to be able to work with elements of $W(\Omega)$. 
\begin{definition} 
\label{def/derivéeNormale}
Let $u \in W(\Omega)$ and $\Delta u \in L_2(\Omega)$.  We say that $u$ has a normal derivative $\psi$ on $\Gamma$, denoted by $\partial_\nu u = \psi$, if $\psi \in L_2(\Gamma)$ and
for all $v \in \mathcal{D}(\R^n)$
\begin{equation}
\label{eq/derivéeNormaleFaible}
\int_{\Omega} (\Delta u)v\dx + \int_{\Omega} \nabla u \cdot \nabla v\dx = \int_\Gamma \psi Tr\,v\dm_d.
\end{equation}
\end{definition}


\begin{remark}
\label{rmq/NormalDerivativeDensity}
Definition~\ref{eq/derivéeNormaleFaible} defines a weak notion of normal derivative of a function in $W(\Omega)$ in the distributional sense, if it exists. If it exists, it is unique. In addition, thanks to the definition of the space $W^D(\Omega)$, functions $v\in \mathcal{D}(\R^n)$, considered on $\Omega$, are dense in $W^D(\Omega)$. Thus, by the density argument, Eq.~(\ref{eq/derivéeNormaleFaible}) holds for all $v \in W^D(\Omega)$ (see~\cite{ARENDT-2015} p.~321).
\end{remark}
Next we define the associated variational formulation for the exterior problem~(\ref{eq/RobCond}):
\begin{definition}
\label{def/weakForm}
Let $\psi \in L_2(\Gamma)$ and $\lambda \geq 0$, we say that $u \in W^D(\Omega)$ is a weak solution to the Robin problem with Dirichlet boundary conditions at  infinity if
\begin{equation}
\label{eq/formFaibleDirichlet}
\forall v \in  W^D(\Omega) \quad \int_\Omega \nabla u \nabla v \dx + \lambda \int_\Gamma \mathrm{Tr} u \mathrm{Tr} v \,\dm_d= \int_\Gamma \psi \mathrm{Tr} v\,\dm_d.
\end{equation}
\end{definition}

The variational formulation~(\ref{eq/formFaibleDirichlet}) is well-posed:
\begin{theorem}\label{ThweakFormWellPosed}
 Let $\Omega$ be an admissible exterior domain 
 with a compact $d$-set boundary $\Gamma$ ($n-2< d<n$, $n\ge 3$).
 For all $\lambda \in [0, \infty[$ and for all $\psi \in L_2(\Gamma)$ there exists a unique weak solution $u\in  W^D(\Omega)$ to the Robin problem with Dirichlet boundary conditions at infinity in the sense of Definition~\ref{def/weakForm}. 
 Moreover, if the operator $B_\lambda$ is defined by $$B_\lambda:\; \psi\in L_2(\Gamma)  \mapsto u\in  W^D(\Omega)$$ with $u$, the solution of Eq.~(\ref{eq/formFaibleDirichlet}), then it satisfies the same properties as the operator $B_\lambda(S)$ introduced in Theorem~\ref{ThWelPTr} for the truncated domains (see points 1--4):
 $B_\lambda$ is a linear compact, positive and monotone operator with $0\le \lambda B_\lambda \mathds{1}_\Gamma\le \mathds{1}_\Omega$ for all $\lambda \in [0, \infty[$.
    \end{theorem}

\textbf{Proof.}
  Thanks to Theorem~\ref{theorem/Ou2},  the trace operator $\mathrm{Tr}$ is continuous from $W^D(\Omega)$ to $L_2(\Gamma)$. Then
 the well-posedness of Eq.~(\ref{eq/formFaibleDirichlet}) and the continuity of $B_\lambda$ follow from the application of the Lax-Milgram theorem in the Hilbert space $W^D(\Omega)$. To prove the compactness of $B_\lambda$, we follow  Ref.~\cite{ARENDT-2015} Proposition 3.9.
 Indeed, let $\lambda \in [0, \infty[$ and $(\psi_k)_{k\in \N}$ be a bounded sequence in $L_2(\Gamma)$. Then there exists $\psi \in L_2(\Gamma)$ such that, up to a sub-sequence, $\psi_k \stackrel{L_2(\Gamma)}{\rightharpoonup} \psi$ for $k\to +\infty$. For all $k \in \N$ we set $u_k = B_\lambda \psi_k$ and $u = B_\lambda \psi$.
 From the continuity of $B_\lambda$ it follows that $u_k \stackrel{W^D(\Omega)}{\rightharpoonup} u$ for $k\to +\infty$. Therefore, $\mathrm{Tr} u_k \stackrel{L_2(\Gamma)}{\rightarrow} \mathrm{Tr} u$ for $k\to +\infty$, since the trace operator $\mathrm{Tr}$ is  compact from $W^D(\Omega)$ to $L_2(\Gamma)$ (see Theorem~\ref{theorem/Ou2} point 4).

Let $k \in \N$, choosing $v = u_k$ in Eq.~(\ref{eq/formFaibleDirichlet}), we obtain
\[
\| u_k \|_{W^D(\Omega)}^2=\int_\Omega |\nabla u_k|^2\dx = \int_\Gamma \psi \mathrm{Tr} u_k\dm_d - \lambda \int_\Gamma |\mathrm{Tr} u_k|^2\dm_d.
\]
Consequently, using Eq.~(\ref{eq/formFaibleDirichlet}) with $v = u$, we have
\begin{equation*}
\lim_{k \rightarrow \infty} \int_\Omega |\nabla u_k|^2\dx = \int_\Gamma \psi \mathrm{Tr} u\, \dm_d- \lambda \int_\Gamma |\mathrm{Tr} u|^2\dm_d
							= \int_\Omega |\nabla u|^2\dx=\| u \|_{W^D(\Omega)}^2.
\end{equation*}
 Hence, $\| u_k \|_{W^D(\Omega)} \rightarrow \| u \|_{W^D(\Omega)}$ for $k\to +\infty$, and  consequently, $B_\lambda$ is compact. 

 The positive and the monotone property of $B_\lambda$ follow respectively from Ref.~\cite{ARENDT-2015} Proposition 3.5 and Proposition 3.7 a). The equality $0\le \lambda B_\lambda \mathds{1}_\Gamma\le \mathds{1}_\Omega$ follows from Ref.~\cite{ARENDT-2015} Proposition 3.6 and Corollary 3.8 b).
$\Box$

Now, let us show that the truncated problem, studied in Subsection~\ref{SecWellPosedTrunc}, independently of the form of the boundary $S$, is an approximation of the exterior problem in $\R^n$ with $n\ge 3$.
We denote by $\Omega_S$ the exterior domain $\Omega$, truncated by the boundary $S$.
In this framework, we also truncate~\cite{ARENDT-2015} the space $W^D(\Omega)$, introducing a subspace
\[
W^D_S(\Omega) := \{ u \in W^D(\Omega) : u|_{\R^n \setminus \Omega_S} = 0 \},
\]
which is closed and, thus, is a Hilbert space for the inner product $(\cdot,\cdot)_{W^D(\Omega)}$.
Since $$H^1_0(\Omega_1) = \{ u|_{\Omega_1} : u \in H^1(\R^n) \mbox{ and } u|_{\R^n \setminus \Omega_1} = 0\},$$ we notice that the map $\Psi: u \in W^D_S(\Omega) \mapsto u|_{\Omega_S}\in \tilde{H}^1(\Omega_S)$ is a bi-continuous bijection. Consequently, problem~(\ref{EqVarFTruncD}) is also well-posed in $W^D_S(\Omega)$ with the same properties described in Theorem~\ref{ThWelPTr}.

In what follows, we will also suppose that the boundary $S$ is far enough from the boundary $\Gamma$. Precisely, we suppose that $\overline{\Omega}_0\subset B_r$ is a domain (all time satisfying the conditions of Theorem~\ref{ThContTrace}), included in a ball of a radius $r_0>0$ (which exists as $\Omega_0$ is bounded), and $\Omega_{S_r}$ with $r\ge r_0$ is such that  $(\R^n\setminus \overline{\Omega}_0)\cap B_r \subset \Omega_{S_r}$ with $\del B_r\cap S_r=\varnothing$. If $r\to +\infty$ the boundaries $S_r$ (for each $r\ge r_0$ the domains $\Omega_{S_r}$ satisfy the conditions of Theorem~\ref{ThContTrace}) will be more and more far from $\Gamma$ and in the limit $r\to +\infty$ the domains $\Omega_{S_r}$ give $\Omega$.
Let us precise the properties of solutions $u\in W^D_S(\Omega)$ for the truncated problem to compare to the solutions on the exterior domain:
\begin{lemma}
\label{lemma/bijection}
Let  $\Omega_0$,  $\Omega$ and $\Omega_S$ (or $\Omega_{S_r}$ for all $r\ge r_0>0$) satisfy conditions of Theorem~\ref{ThCompTrExtdom} in $\R^n$ with $n\ge 3$.
Let $B_\lambda(S): \psi \in L_2(\Gamma)\mapsto u\in W^D_S(\Omega)$ be the operator for the truncated problem and $B_\lambda: \psi \in L_2(\Gamma)\mapsto u\in W^D(\Omega)$ be the operator for the exterior problem.

Then for all $\lambda\in [0,\infty[$ and $\psi\in L_2(\Gamma)$, if $\psi\ge 0$ in $L_2(\Gamma)$  and $r_2\ge r_1\ge r_0$, then $$0\le u_{S_{r_1}}=B_\lambda(S_{r_1})\psi\le u_{S_{r_2}}= B_\lambda(S_{r_2})\psi\le B_\lambda \psi=u.$$
\end{lemma}
\textbf{Proof.}
The proof follows the analogous proof as in Ref.~\cite{ARENDT-2015} Propositions 3.5 and 3.6 (see also~\cite{ARENDT-2015} Proposition 4.4).
$\Box$
We can now state the approximation result, ensuring that a solution in any admissible truncated domain, even with a fractal boundary, but which is sufficiently far from $\Gamma$ is an approximation of the solution of the exterior problem:
\begin{theorem}
\label{theorem/convOpRobin}
Let $\lambda \in [0, \infty)$, $\psi\in L_2(\Gamma)$ and $(S_m)_{m\in \N}$ be a fixed sequence of the boundaries of the truncated domains $\Omega_{S_m}$ in $\R^n$ ($n\ge 3$), satisfying for all $m\in \N$ the conditions of Theorem~\ref{ThCompTrExtdom} and such that $(\Omega_{S_m}\cup \overline{\Omega}_0)\supset B_m \supset \overline{\Omega}_0$. Let $u_{S_m}=B_\lambda(S_m)\psi$ and $u= B_\lambda \psi$. Then for all $\eps>0$ there exits  $m_0(\eps)>0$,  independent on  the chosen sequence of the boundaries  $(S_m)$, such that  $$\forall m\ge m_0
\quad \|u_{S_m}-u\|_{W^D(\Omega)}<\eps.$$ Equivalently, for all described sequences $(S_m)_{m\in \N}$,  it holds
 $$\|B_\lambda(S_m) - B_\lambda\|_{\mathcal{L}(L_2(\Gamma), W^D(\Omega))}\to 0\quad \hbox{as} \quad m\to +\infty.$$ 
\end{theorem}
\textbf{Proof.}It is a simple generation using our previous results of  Theorem 4.3~\cite{ARENDT-2015}. 
$\Box$

 \section{Spectral properties of the Poincar\'e-Steklov operator defined  by the interior and by the exterior problems}\label{SecDtoNReg}
 The Poincar\'e-Steklov operator, also named the Dirichlet-to-Neumann operator, was originally introduced by V.A. Steklov and usually defined by a map
 $$A: u|_\Gamma \mapsto \frac{\del u }{\del \nu}|_\Gamma$$ for a solution $u$ of the elliptic Dirichlet problem:
 $-\Delta u=0$ in a domain $\Omega$ and $u|_\Gamma=f$ ($\del \Omega=\Gamma$).
 
 It is well-known that if $\Omega$ is a bounded domain with $C^\infty$-regular boundary (a regular manifolds with boundary), then the operator $A: C^\infty(\Gamma)\to C^\infty(\Gamma)$ is an elliptic self-adjoint pseudo-differential operator of the first order (see~\cite{TAYLOR-1996} $\S$11 and 12 of Chapter 7) with a discrete spectrum
 $$0=\lambda_0<\lambda_1\le \lambda_2\le \ldots, \quad \hbox{with } \lambda_k\to +\infty \; k\to +\infty.$$
 If $A$ is considered as an operator $H^1(\Gamma)\to L_2(\Gamma)$, then  its eigenfunctions form a basis in $L_2(\Gamma)$.
 For any Lipschitz boundary $\Gamma$ of a bounded domain $\Omega$, the Dirichlet-to-Neumann operator $$A: H^\frac{1}{2}(\Gamma)\to H^{-\frac{1}{2}}(\Gamma)$$ is well-defined and it is a linear continuous self-adjoint operator. 
Thanks to~\cite{ARENDT-2012}, we also know that 
 the Dirichlet-to-Neumann operator $A$ has compact resolvent, and hence discrete spectrum, as long as the trace operator  $\mathrm{Tr}: H^1(\Omega) \to L_2(\Gamma)$ is compact (see also~\cite{ARENDT-2011} and~\cite{TRIEBEL-1997} for abstract definition of the elliptic operators on a $d$-set). Thus, thanks to Theorem~\ref{ThCompTrExtdom}, the property of the compact resolvent also holds for an admissible $n$-set $\Omega$ with a compact $d$-set boundary $\Gamma$. We will discuss it in details in the next section. From~\cite{BEHRNDT-2015}, we also have that  $\operatorname{Ker}{A}\ne \{0\}$, since $0$ is the eigenvalue of the Neumann eigenvalue problem for the Laplacian.
 For the Weil asymptotic formulas for the distribution of the eigenvalues of the Dirichlet-to-Neumann operator there are results for bounded smooth compact Riemannian manifolds with $C^\infty$ boundaries~\cite{GIROUARD-2014}, for polygons~\cite{ARXIV-GIROUARD-2014} and more general class of plane domains~\cite{GIROUARD-2015} and also for a bounded domain with a fractal boundary~\cite{PINASCO-2005}.
%
%
%

 
 In the aim to relate these spectral results, obtained for the Dirichlet-to-Neumann operator for a bounded domain, with the case of the exterior domain, we prove the following theorem:
 \begin{theorem}\label{ThSPECTR}
   Let $\Omega_0\subset \R^n$ ($n\ge 2$) be an admissible bounded domain with a compact boundary $\Gamma$ such that its complement in $\R^n$ $\Omega=\R^n\setminus \overline{\Omega}_0$ be also an admissible domain with the same boundary $\Gamma$, satisfying conditions of Theorem~\ref{ThCompTrExtdom}. 
   Then the Dirichlet-to-Neumann operators $$A^{int}: L_2(\Gamma) \to L_2(\Gamma), \hbox{ associated with the Laplacian on } \Omega_0,$$  
   and $$A^{ext}: L_2(\Gamma) \to L_2(\Gamma),\hbox{ associated with the Laplacian on } \Omega,$$  
  are self-adjoint positive operators with  compact resolvents and  discrete positive spectra. Let us denote the sets of all eigenvalues of $A^{int}$ and $A^{ext}$ respectively by $\sigma^{int}$ and $\sigma^{ext}$, which are subsets of $\R^+$. 
   
   If $\Omega_{S}$ is an admissible truncated domain to $\Omega_0$, then the associated Dirichlet-to-Neumann operator $$A(S): L_2(\Gamma) \to L_2(\Gamma)$$ for all $n\ge2$  is  self-adjoint positive operator with a compact resolvent and a discrete spectrum. The point spectrum, or the set of all eigenvalues of $A(S)$, is strictely positive:
   $\sigma_S\subset \R^+_*$ ($i.e.$ $A(S)$ is injective with the compact inverse operator $A^{-1}(S)$). 
   
   In addition, let  $\mu_k(r)\in \sigma_S(r)$, where $\sigma_S(r)\subset \R^+_*$ is the point spectrum of the Dirichlet-to-Neumann operator $A(S_r)$, associated with an admissible truncated domain $\Omega_{S_r}$, such that $(\Omega_{S_r}\cup \overline{\Omega}_0)\supset B_r \supset \overline{\Omega}_0.$                
             For $n=2$,  if $\mu_0(r)=\min_{k\in \N} \mu_k(r)$, then, independently on the form of $S_r$, 
             \begin{equation}\label{EqT1}
                                                                      \mu_0(r)\to 0\quad \hbox{for} \quad r\to +\infty.
                                                                     \end{equation}

             For $n\ge 3$, \begin{equation}\label{EqT2} 
                            \|A^{-1}(S_r)-A^{-1}\|_{\mathcal{L}(L_2(\Gamma))}\to 0\quad \hbox{for} \quad r\to +\infty
                           \end{equation}
 independently on the form of $S_r$.
   
   Moreover, 
all non-zero eigenvalue of $A^{int}$ 
is also  an eigenvalue of $A^{ext}$ and converse. Hence the eigenfunctions of $A^{int}$ and $A^{ext}$ form the same basis in $L_2(\Gamma)$.

More precisely, 
it holds
\begin{itemize}
                \item For $n=2$   $$\sigma^{int}=\sigma^{ext}\subset \R^+\quad  \hbox{and } \quad 0\in \sigma^{ext}.$$ 
                \item For $n\ge 3$  $$\sigma^{int}=\{0\}\cup \sigma^{ext}\quad \hbox{ with }\quad \sigma^{ext}\subset ]0,+\infty[,$$ $i.e.$ the Dirichlet-to-Neumann operator of the exterior problem, also as of the truncated problem, is an injective operator with the compact inverse. 
                
                                                                                                                 \end{itemize}
                                                                                                                 \end{theorem}
                                                                                                                 
%
 To prove Theorem~\ref{ThSPECTR}, we need to define the Dirichlet-to-Neumann operator on a $d$-set $\Gamma$ in $L_2(\Gamma)$. Hence, we firstly do it in Section~\ref{SecDtoNFrac} and then give the proof in Section~\ref{SecFinal}.
 


  \section{Poincar\'e-Steklov operator on a $d$-set}\label{SecDtoNFrac}
 \subsection{For a bounded domain}\label{SsecDtoNBound}
 Let $\Omega_0$ be a bounded admissible domain with a $d$-set boundary $\Gamma$ ($n-2< d<n$, $n\ge 2$). 
 Knowing the well-posedness results for the Dirichlet problem (Theorem 7~\cite{JONSSON-1997}) and the definition of the normal derivative by the Green formula~(\ref{FracGreen}), thanks to~\cite{JONSSON-1995}, we notice that the general setting of~\cite{BEHRNDT-2015} p. 5904 for Lipschitz domains (see also~\cite{MCLEAN-2000} Theorem 4.10) still holds in the  case of a $d$-set boundary, by replacing 
 $H^\frac{1}{2}(\Gamma)$ by $B^{2,2}_\beta(\Gamma)$ with $\beta=1-\frac{n-d}{2}>0$ and $H^{-\frac{1}{2}}(\Gamma)$ by $B^{2,2}_{-\beta}(\Gamma)$.
  Precisely, we have
  that for all $\lambda\in \C$ the Dirichlet problem
  \begin{equation}\label{}
   -\Delta u=\lambda u, \quad u|_\Gamma=\phi
  \end{equation}
 is solvable if $\phi\in B^{2,2}_\beta(\Gamma)$  satisfies
  $$\langle\del_\nu \psi,\phi\rangle_{B^{2,2}_{-\beta}(\Gamma)\times B^{2,2}_\beta(\Gamma)}=0$$ for all solutions $\psi\in H^1(\Omega_0)$ of the corresponding homogeneous problem
  \begin{equation}\label{EqDual}
   -\Delta \psi=\lambda \psi, \quad \psi|_\Gamma=0.
  \end{equation}
We are especially interesting in the case $\lambda=0$.
Thus, we directly conclude that  problem~(\ref{EqDual}) has only the trivial solution $\psi=0$ ($\lambda=0$ is not an eigenvalue of the Dirichlet Laplacian), and consequently the Poincar\'e-Steklov  operator
$A:B^{2,2}_\beta(\Gamma)\to B^{2,2}_{-\beta}(\Gamma)$ mapping $u|_\Gamma$ to $\del_\nu u|_\Gamma$ is well-defined on $B^{2,2}_\beta(\Gamma)$.

On the other hand, as it was done in~\cite{ARENDT-2011} for bounded domains with $(n-1)$-dimensional boundaries, it is also possible to consider $A$ as operator from $L_2(\Gamma)$ to $L_2(\Gamma)$, if we consider the trace map $\mathrm{Tr}: H^1(\Omega_0)\to L_2(\Gamma)$ (note that $B^{2,2}_\beta(\Gamma)\subset L_2(\Gamma)$) and update the definition of the normal derivative by analogy with Definition~\ref{def/derivéeNormale}:
\begin{definition}\label{DefDerNormBound}
 Let $u \in H^1(\Omega_0)$ and $\Delta u \in L_2(\Omega_0)$. If there exists $\psi \in L_2(\Gamma)$ such that for all $v \in H^1(\Omega_0)$ it holds Eq.~(\ref{eq/derivéeNormaleFaible}), then $\psi$ is called  a  $L_2$-normal derivative of $u$, denoted by   $\partial_\nu u = \psi$. 
\end{definition}
Definition~\ref{DefDerNormBound} restricts the normal derivative of $u$, which is naturally in $B^{2,2}_{-\beta}(\Gamma)$, to a consideration of only the normal derivative from its dense subspace. Thus, the $L_2$-normal derivative can does not exist, but if it exits, it is unique.

Therefore, to define the Dirichlet-to-Neumann operator on $L_2(\Gamma)$, we use the following Theorem from~\cite{ARENDT-2011} (see Theorem 3.4)
\begin{theorem}\label{thmOpaj}
Let $D(a)$ be a real vector space and let $a : D(a) \times D(a) \rightarrow \R $ be bilinear symmetric such that $a(u, u) \geq 0$ for all $u \in D(a)$. Let $H$ be a real Hilbert space and let $T : D(a) \rightarrow H$ be linear operator with dense image. Then there exists a positive and self-adjoint operator $A$ on $H$ such that for all $\phi,\, \psi \in H$, one has $\phi \in D(A)$ and $A\phi = \psi$ if and only if there exists a sequence $(u_k)_{k \in \N }$ in $D(a)$ such that:
\begin{enumerate}
\item $\lim_{k, m \rightarrow \infty} a(u_k - u_m,u_k - u_m) = 0$,
\item $\lim_{k \rightarrow \infty} T(u_k) = \phi$ in $H$,
\item for all $v \in D(a)$  $\lim_{k \rightarrow \infty} a(u_k, v) = (\psi, T(v))_H$. 
\end{enumerate}
 The operator $A$ is called the operator associated with $(a, T)$.
\end{theorem}

 Consequently we state
 \begin{theorem}\label{ThDtoNInt}
  Let $\Omega_0$ be a bounded admissible domain with a compact $d$-set boundary $\Gamma$ ($n-2< d<n$, $n\ge 2$). Then for  $\beta=1-\frac{n-d}{2}>0$ the Poincar\'e-Steklov operator $$A:B^{2,2}_\beta(\Gamma)\to B^{2,2}_{-\beta}(\Gamma)$$ mapping $u|_\Gamma$ to $\del_\nu u|_\Gamma$ is linear bounded 
  self-adjoint operator with $\operatorname{Ker} A\ne \{0\}$. 
  In addition, the Poincar\'e-Steklov operator $A$, considered  from $L_2(\Gamma)$ to $L_2(\Gamma)$, is self-adjoint positive operator with a compact resolvent.  Therefore, there exists a discrete spectrum of $A$ with eigenvalues
  $$0=\mu_0<\mu_1\le \mu_2\le \ldots, \quad \hbox{with } \mu_k\to +\infty \; k\to +\infty$$
  and the corresponding eigenfunctions form an orthonormal basis in $L_2(\Gamma)$.
 \end{theorem}
\textbf{Proof.}
We have already noticed that the domain of $A$ is exactly $B^{2,2}_\beta(\Gamma)$. As $0$ is an eigenvalue of the Neumann Laplacian, $\operatorname{Ker} A\ne \{0\}$.
From the following Green formula for $u,$ $v\in H^1(\Omega_0)$ with $\Delta u,$ $\Delta v\in L_2(\Omega_0)$
\begin{multline*}
 \int_{\Omega_0} \Delta uv \dx-\int_{\Omega_0} u\Delta v \dx\\=\langle\frac{\del u}{\del \nu}, \mathrm{Tr} v \rangle_{B^{2,2}_{-\beta}(\Gamma),B^{2,2}_\beta(\Gamma)}- \langle \mathrm{Tr} u, \frac{\del v}{\del \nu}\rangle_{B^{2,2}_\beta(\Gamma),B^{2,2}_{-\beta}(\Gamma)},
\end{multline*}
we directly find that for all $u,$ $v\in B^{2,2}_\beta(\Gamma)$
$$\langle A u,  v \rangle_{B^{2,2}_{-\beta}(\Gamma),B^{2,2}_\beta(\Gamma)}= \langle  u, Av\rangle_{B^{2,2}_\beta(\Gamma),B^{2,2}_{-\beta}(\Gamma)},$$
 $i.e.$ the operator $A$ is self-adjoint and closed. Since  $B^{2,2}_{-\beta}(\Gamma)$ is a Banach space, by the closed graph Theorem, $A$ is continuous.
 
 To define $A$ as an operator on $L_2(\Gamma)$ we use~\cite{ARENDT-2011,ARENDT-2012,ARENDT-2012-1}. As $\Omega_0$ is such that the trace operator $\mathrm{Tr}$ is compact from $H^1(\Omega_0)$ to $L_2(\Gamma)$, then the embedding of its image $\mathrm{Tr}(H^1(\Omega_0))=B^{2,2}_\beta(\Gamma)$ into $L_2(\Gamma)$ is compact. Now, as it was noticed in~\cite{ARENDT-2011}, the space $\{v|_{\Gamma} : v \in \mathcal{D}(\R ^n)\}$ is dense in $C(\Gamma)$ by the Stone-Weierstrass theorem for the uniform norm  and, therefore, it is also dense in $L_2(\Gamma)$, since we endowed $\Gamma$ with the $d$-dimensional Hausdorff measure which is Borel regular. Hence, $B^{2,2}_\beta(\Gamma)$ is dense in $L_2(\Gamma)$.  It allows us to apply Theorem 2.2 and follow Section 4.4 of Ref.~\cite{ARENDT-2012-1}.
 
 Using the results of Section~\ref{SecFuncAnT}, we follow the proof of Wallin~\cite{WALLIN-1991}, Theorem 3, to obtain that for all bounded admissible domains
 $$\operatorname{Ker} \mathrm{Tr}=H_0^1(\Omega_0).$$
  Thanks to Lemma 2.2~\cite{ARENDT-2012}, $$H^1(\Omega_0)=H^1_0(\Omega_0)\oplus H$$ with $$H=\{u\in H^1(\Omega_0)|\; \Delta u=0 \hbox{ weakly}\}.$$ Hence, $\mathrm{Tr}(H)=B^{2,2}_\beta(\Gamma)$ and $\mathrm{Tr}:H\to B^{2,2}_\beta(\Gamma)$ is a linear bijection.
 Therefore, the bilinear map $a_0: B^{2,2}_\beta(\Gamma) \times B^{2,2}_\beta(\Gamma) \to \R$, given by
 \begin{equation}\label{EqBilFormn2}
  a_0(\phi,\psi)=\int_{\Omega_0} \nabla u \nabla v\dx \quad \hbox{ for } u,v \in H \quad \mathrm{Tr}u=\phi,\quad \mathrm{Tr}v=\psi,
 \end{equation}
is symmetric, continuous and elliptic~\cite{ARENDT-2007} (see Proposition 3.3, based on the compactness of the embedding $H\subset\subset L_2(\Omega_0)$ (as $H$ is a closed subspace of $H^1(\Omega)$ and $H^1(\Omega)\subset\subset L_2(\Omega_0)$, this implies $H\subset\subset L_2(\Omega_0)$)  and on the injective property of the trace from $H$ to $L_2(\Gamma)$):
$$\exists \omega\ge 0\hbox{ such that} \quad \forall u\in H \quad a_0(\mathrm{Tr}u,\mathrm{Tr}u)+\omega\int_\Gamma |u|^2\dm_d\ge \frac{1}{2}\|u\|^2_{H^1(\Omega_0)}.$$
If the operator $N: L_2(\Gamma)\to L_2(\Gamma)$ is the operator  associated with $a_0$, then it is the Dirichlet-to-Neumann operator $A$ on $L_2(\Gamma)$, $i.e.$ $A\phi=\del_\nu u$ in $L_2(\Gamma)$ with 
\begin{multline*}
 \mathcal{D}(A)=\{\phi\in L_2(\Gamma)|\; \exists u\in H^1(\Omega_0) \hbox{ such that } \\\mathrm{Tr}u=\phi, \quad \Delta u=0 \hbox{ and } \exists \del_\nu u \in L_2(\Gamma)\}.
\end{multline*}
 Moreover, we have  that  for all $\phi\in L_2(\Gamma)$,  \begin{multline*}
                                                           \phi \in D(A) \hbox{ and there exists  an element } \psi=A\phi \hbox{ of } L_2(\Gamma)\quad \Longleftrightarrow\\
                                                           \exists u\in H^1(\Omega_0)\hbox{ such that }\mathrm{Tr}u=\phi \hbox{ and }\forall v\in H^1(\Omega_0)\quad \int_{\Omega_0}\nabla u\nabla v\dx=\int_\Gamma \psi \mathrm{Tr}v\dm_d.
                                                          \end{multline*}

 On the other hand, we also can directly use Theorem~3.3 in Ref.~\cite{ARENDT-2011}, by applying Theorem~\ref{thmOpaj}.
 Let now $D(a)=H^1(\Omega_0)\cap C(\overline{\Omega}_0)$, which is dense in $H^1(\Omega_0)$ (see the discussion of Ref.~\cite{ARENDT-2011}). Then $\mathrm{Tr}(D(a))$ is dense in $L_2(\Gamma)$.
 Therefore, taking in Theorem~\ref{thmOpaj} $$a(u,v)=\int_{\Omega_0} \nabla u\nabla v\dx: D(a)\times D(a)\to \R, \; H=L_2(\Gamma) \hbox{ and } T=\mathrm{Tr}: D(a)\to L_2(\Gamma),$$ as $\mathrm{Tr}$ is compact, we conclude that the operator associated to $(a,\mathrm{Tr})$ is the Dirichlet-to-Neumann operator $A$, positive and self-adjoint in $L_2(\Gamma)$ (see the proof of Theorem~3.3 in Ref.~\cite{ARENDT-2011}). 
  Since the compactness of the trace implies that $A$ has a compact resolvent,  it is sufficient to apply the Hilbert-Schmidt Theorem to finish the proof.
$\Box$

       \subsection{For an exterior and truncated domains}

       In this subsection we generalize~\cite{ARENDT-2015} and introduce the Dirichlet-to-Neumann operator $A$ on $L_2(\Gamma)$ with respect to the exterior domain $\Omega\subset \R^n$ and $A(S)$ with respect to a truncated domain for $n\ge 2$ in the framework of $d$-sets. 

\begin{definition}\textbf{(Dirichlet-to-Neumann operator for an exterior domain $n\ge 3$)}
\label{def/opDtN}
 Let $\Omega\subset \R^n$, $n\ge 3$, be an  admissible exterior domain, satisfying the conditions of Theorem~\ref{ThCompTrExtdom}. 
 The operator $A: L_2(\Gamma)\to L_2(\Gamma)$,  associated with the bilinear form
$a^D : W^D(\Omega) \times W^D(\Omega) \rightarrow \R$ given by 
$$a^D(u, v) = \int_\Omega \nabla u \nabla v\dx=\langle u,v\rangle_{W^D(\Omega)},$$
and the trace operator $\mathrm{Tr}: W^D(\Omega)\to L_2(\Gamma)$, is called the Dirichlet-to-Neumann operator with the Dirichlet boundary condition at infinity. 
\end{definition}
\begin{remark}\label{RemLaplExtN2}
 Theorem~\ref{thmOpaj} does not require to $D(a)$ the completeness, $i.e.$ $a(\cdot,\cdot)$ can be equivalent to a semi-norm on $D(a)$, what is the case of $W^D(\Omega)$ with $a(u,u)=\int_\Omega |\nabla u|^2 \dx$ for $n=2$. Therefore, it allows us to define the Dirichlet-to-Neumann operator $A$ of the exterior problem in $\R^2$, which can be understood as the limit case for $r\to +\infty$ of the problem for a truncated domain  well-posed in $\tilde{H}^1(\Omega_{S_r})$. 
 In the case of $W^D(\Omega)$ in $\R^n$ with $n\ge 3$, we have that $D(a)=W^D(\Omega)$ is the Hilbert space corresponding to the inner product $a(\cdot,\cdot)$.
\end{remark}
Let us notice that the trace on the boundary $\Gamma$ satisfies $$\mathrm{Tr}(\mathcal{D}(\R^n))\subset \mathrm{Tr}(W^D(\Omega))\subset L_2(\Gamma)$$ and, since $\mathrm{Tr}(\mathcal{D}(\R^n))$ is dense in $L_2(\Gamma)$, $\mathrm{Tr}(W^D(\Omega))$ is dense in $L_2(\Gamma)$. 
In addition, $a^D$ is $\mathrm{Tr}$-elliptic thanks to Point 2 of Theorem~\ref{theorem/Ou2}, 
$i.e.$
there exists $\alpha\in \R$ and $\delta>0$ such that
$$\forall u\in W^D(\Omega)\quad \int_\Omega |\nabla u|^2\dx+\alpha\int_\Gamma |\mathrm{Tr}u|^2\dm_d\ge \delta \int_\Omega |\nabla u|^2\dx.$$
Thus, for $n\ge 3$ we can also apply Theorem 2.2 and follow Section 4.4 of Ref.~\cite{ARENDT-2012-1}.

For the two-dimensional case, we define $A$ associated to the bilinear form $a_0$ from Eq.~(\ref{EqBilFormn2}), initially given for the interior case:
\begin{theorem}\label{PropDtNn2}\textbf{(Dirichlet-to-Neumann operator for an exterior domain $n=2$)}
 Let $\Omega\subset \R^2$ be an admissible exterior domain, satisfying the conditions of Theorem~\ref{ThCompTrExtdom}. 
 The operator $A: L_2(\Gamma)\to L_2(\Gamma)$,  associated with the bilinear form $a_0$, defined in Eq.~(\ref{EqBilFormn2}), is the Dirichlet-to-Neumann operator with the Dirichlet boundary condition at infinity in the sense that for all $\phi\in L_2(\Gamma)$, 
 \begin{multline*}
            \phi \in D(A)  \hbox{ and there exists  an element }        \psi=A\phi\in    L_2(\Gamma) \quad \Longleftrightarrow \\            
                                      \exists u\in H^1(\Omega) \hbox{ such that }\mathrm{Tr}u=\phi \hbox{ and }\forall v\in H^1(\Omega)\quad \int_\Omega\nabla u\nabla v\dx=\int_\Gamma \psi \mathrm{Tr}v\dm_d.                                                                                                                                                                                                                                \end{multline*}
       Therefore, the properties of $A$ are the same as for the bounded domain case in Theorem~\ref{ThDtoNInt}: the Poincar\'e-Steklov operator $A$ is self-adjoint positive operator with a compact resolvent, and  a discrete spectrum containing positive eigenvalues
  $$0=\mu_0<\mu_1\le \mu_2\le \ldots, \quad \hbox{with } \mu_k\to +\infty \quad \hbox{for } k\to +\infty.$$
  The corresponding eigenfunctions form an orthonormal basis in $L_2(\Gamma)$.                              
\end{theorem}
\textbf{Proof.}
 We use  that $H^1(\R^2)=H^1_0(\R^2)$ and that the compactness of the embedding $H=\{u\in H^1(\Omega)|\; \Delta u=0\hbox{ weakly}\}\subset L_2(\Omega)$ and the injective property of the trace from $H$ to $L_2(\Gamma)$ still hold for the exterior case. In addition $0$ is not an eigenvalue of the Dirichlet Laplacian on $\Omega$. Thus we can follow the proof of Lemma 3.2 and Proposition 3.3  in~\cite{ARENDT-2007}, given for a Lipschitz bounded domain. The spectral properties of $A$ are deduced from the analogous properties proved in Theorem~\ref{ThDtoNInt}.
$\Box$

The following proposition legitimates Definition~\ref{def/opDtN} in the framework of Theorem~\ref{thmOpaj} for $n\ge 3$:
\begin{proposition}
\label{prop/opDtN}
Let $\Omega\subset \R^n$, $n\ge 3$, be an admissible exterior domain, satisfying the conditions of Theorem~\ref{ThCompTrExtdom}, and 
let $\phi$, $\psi \in L_2(\Gamma)$. Then 

$\phi \in D(A)$ and $A\phi = \psi$ if and only if there exists a function $u \in W^D(\Omega)$ such that  $\mathrm{Tr} u = \phi$, $\Delta u = 0$ weakly and $\partial_\nu u = \psi$ in the sense of Definition~\ref{def/derivéeNormale}.
\end{proposition}
\textbf{Proof.}
 Let $\phi$, $\psi \in L_2(\Gamma)$ such that $\phi \in D(A)$ and $A\phi = \psi$. Then, according to Theorem~\ref{thmOpaj}, there exists a sequence $(u_k)_{k \in \N }$ in $W^D(\Omega)$ such that
\begin{enumerate}
\item $\lim_{k, m \rightarrow \infty} \int_\Omega |\nabla (u_k - u_k)|^2 \dx= 0$,
\item $\lim_{k \rightarrow \infty} \mathrm{Tr} u_k = \phi$,
\item $\lim_{k \rightarrow \infty} \int_\Omega \nabla u_k \nabla v \dx= \int_\Gamma \psi \mathrm{Tr}v\dm_d$ for all $v \in W^D(\Omega)$.
\end{enumerate}
Form Item~$1$ it follows  that $(u_k)_{k \in \N }$ is a Cauchy sequence in $W^D(\Omega)$. Therefore, by  the completeness of  $W^D(\Omega)$ (thanks to $n\ge 3$),  there exists $u \in W^D(\Omega)$ such that $u_k \to u$ in $W^D(\Omega)$.   Moreover, since $\mathrm{Tr}: W^D(\Omega)\to L_2(\Gamma)$  is continuous by Point~4 of Theorem~\ref{theorem/Ou2}, $\mathrm{Tr} u = \phi$, according to Item~2.
From Item~3 we deduce that for all $v \in W^D(\Omega)$
\begin{equation}
\label{eq/preuve/propOpDtN}
\int_\Omega \nabla u \nabla v\dx = \int_\Gamma \psi \mathrm{Tr}v\dm_d,
\end{equation}
 and hence, in particular for all $v \in \mathcal{D}(\Omega)$. Therefore $\Delta u = 0$. This with Eq.~(\ref{eq/preuve/propOpDtN}) yields that $u$ has a normal derivative in $L_2(\Gamma)$ and $\partial_\nu u = \psi$.

Conversely, let $\phi$, $\psi \in L_2(\Gamma)$  be such that there exists a function $u \in W^D(\Omega)$, so that  $\mathrm{Tr} u = \phi$, $\Delta u = 0$, $\partial_\nu u = \psi$. According to the definition of normal derivatives (see Definition~\ref{def/derivéeNormale} and Remark~\ref{rmq/NormalDerivativeDensity}), since  $\Delta u = 0$, we have for all $v \in W^D(\Omega)$:
\[
\int_\Omega \nabla u \nabla v\dx = \int_\Gamma \psi \mathrm{Tr}v\dm_d.
\]
 Therefore, for $n\ge 3$ we can apply  Theorem~\ref{thmOpaj} to the sequence, defined by $u_k = u$ for all $k \in \N $, and the result follows.
$\Box$
Let us notice that the Dirichlet-to-Neumann operator $A(S)$ for a domain, truncated by a $d_S$-set $S$ ($n-1\le d_S<n$), can be defined absolutely in the same way as the operator $A$ for the exterior domains if we replace $W^D(\Omega)$ by $\tilde{H}^1(\Omega_S)$ or, equivalently, by $W^D_S(\Omega)$.

Consequently, for  exterior and truncated domains we have
\begin{theorem}
\label{ThDtNcompResolventeExt}
Let $\Omega$ be an admissible exterior domain in $\R^n$ with $n\ge 3$  and $\Omega_S$ be an admissible truncated domain in $\R^n$ with $n\ge 2$ satisfying conditions of Theorem~\ref{ThCompTrExtdom} and $\lambda \in [0, \infty[$. Then the Dirichlet-to-Neumann operator with the Dirichlet boundary condition at infinity $A: L_2(\Gamma)\to L_2(\Gamma)$ (see Definition~\ref{def/opDtN})  and the Dirichlet-to-Neumann operator $A(S)$ of the truncated domain are positive self-adjoint operators  with a compact resolvent $$\forall \lambda\in[0,+\infty[\quad   (\lambda I + A)^{-1}= \mathrm{Tr}_\Gamma\circ B_\lambda, \quad (\lambda I + A(S))^{-1}= \mathrm{Tr}_\Gamma\circ B_\lambda(S) $$ where $B_\lambda: \; \psi\in L_2(\Gamma)  \mapsto u\in  W^D(\Omega)$ with $u$, the solution of Eq.~(\ref{eq/formFaibleDirichlet}), and $B_\lambda(S): \; \psi\in L_2(\Gamma)  \mapsto u\in  W^D_S(\Omega)$ with $u$, the solution of Eq.~(\ref{EqVarFTruncD})   are defined in Theorem~\ref{ThweakFormWellPosed} and Theorem~\ref{ThWelPTr} respectively. Moreover, $\operatorname{Ker} A= \operatorname{Ker} A(S)= \{0\}$ and for $n\ge 3$, independently on a $d$-set $S_r$ ($(\R^n\setminus \overline{\Omega}_0)\cap B_r \subset \Omega_{S_r}$ with $\del B_r\cap S_r=\varnothing$),
$$\forall \lambda\in[0,+\infty[ \quad\|(\lambda I+A(S_r))^{-1}-(\lambda I+A)^{-1}\|_{\mathcal{L}(L_2(\Gamma))}\to 0 \hbox{ as } r\to +\infty.$$

Therefore, the spectra of $A$ ($n\ge 3$) and $A(S)$ ($n\ge 2$) are discrete with all eigenvalues $(\mu_k)_{k\in\N}$ (precisely, $(\mu_k(A))_{k\in \N}$ of $A$ and $(\mu_k(A(S)))_{k\in \N}$ of $A(S)$) strictly positive 
$$0<\mu_0<\mu_1\le \mu_2\le \ldots, \quad \hbox{with } \mu_k\to +\infty \quad \hbox{for }  k\to +\infty,$$
and the corresponding eigenfunctions form  orthonormal basis of $L_2(\Gamma)$.
\end{theorem}
\textbf{Proof.}
 The compactness of the resolvents $(\lambda I + A)^{-1}$ and $(\lambda I + A(S))^{-1}$ directly follows from the compactness properties of the operators $\mathrm{Tr}_\Gamma$, $B_\lambda$, $B_\lambda(S)$. Using the previous results and the Hilbert-Schmidt Theorem for self-adjoint compact operators on a Hilbert space, we finish the proof. 
$\Box$

\section{Proof of Theorem~\ref{ThSPECTR} and final remarks}\label{SecFinal}

Now, we can prove Theorem~\ref{ThSPECTR}:
\textbf{Proof.}
  Actually, Theorems~\ref{ThDtoNInt},~\ref{PropDtNn2} and~\ref{ThDtNcompResolventeExt} implies that the operators $A^{int}$, $A^{ext}$ and $A(S_r)$ have compacts resolvents and discrete positive spectra. As previously, by $\sigma^{int}$, $\sigma^{ext}$ and $\sigma_S(r)$ are denoted  the sets of all eigenvalues of $A^{int}$, $A^{ext}$ and $A(S_r)$ respectively.
  With these notations, for all $n\ge 2$ the point $0\notin \sigma_S(r)$ (by Theorem~\ref{ThDtNcompResolventeExt}), but $0\in \sigma^{int}$ (by Theorem~\ref{ThDtoNInt}). Thanks to Theorem~\ref{PropDtNn2}, for $n=2$ the point $0\in   \sigma^{ext}$,  and,  thanks to Theorem~\ref{ThDtNcompResolventeExt}, for $n\ge 3$ the point $0\notin \sigma^{ext}$. 
  The approximation result for the resolvents of the exterior and truncated domains in Theorems~\ref{ThDtNcompResolventeExt} (for $\lambda=0$)  gives
Eqs.~(\ref{EqT2}). 

  
  Thus, we need to prove that all non-zero eigenvalues of $A^{int}$ 
are also   eigenvalues of $A^{ext}$ and converse.  Grebenkov~\cite{GREBENKOV-2004} (pp. 129-132 and 134) have shown it by the  explicit calculus of the interior and exterior spectra of the Dirichlet-to-Neumann operators  for a ball.

If $\Gamma$ is regular, it is sufficient to apply a conform map to project $\Gamma$ to a sphere and, hence, to obtain the same result (for the conformal map technics see~\cite{GIROUARD-2014} the proof of Theorem 1.4, but also~\cite{GIROUARD-2015} and~\cite{BANJAI-2007}). 
For the general case of a $d$-set $\Gamma$, it is more natural to use  given in the previous Section definitions of the Dirichlet-to-Neumann operators.

Let $n\ge 3$. If $\mu>0$ is an eigenvalue of $A^{ext}$, corresponding to an eigenfunction $\phi\in L_2(\Gamma)$, then, according to Proposition~\ref{prop/opDtN},
\begin{multline}
 \phi \in D(A)\hbox{ and } A\phi =\mu \phi\hbox{ if and only if } \exists u \in W^D(\Omega) \hbox{ such that}\\
 \mathrm{Tr} u = \phi,\quad  \Delta u = 0\hbox{ and } \partial_\nu u = \mu \phi, \; i.e.\\
 \forall v\in \mathcal{D}(\R^n) \quad \int_\Omega\nabla u\nabla v\dx=\int_{\Gamma}\mu \phi \mathrm{Tr}v\dm_d. 
\end{multline}

The trace on $v$ on $\Gamma$ can be also considered for a function $v\in H^1(\Omega_0)$, and, by the same way, $\phi\in L_2(\Gamma)$ can be also interpreted as the trace of $w\in H^1(\Omega_0)$. Thus,  $\mu \phi\in L_2(\Gamma)$ is a normal derivative of $w\in H^1(\Omega_0)$ if and only if
$$\forall v\in H^1(\Omega_0) \quad \int_{\Omega_0}\nabla w\nabla v\dx=\int_{\Gamma}\mu \phi \mathrm{Tr}v\dm_d \hbox{ and } \Delta w=0 \hbox{ weakly in } \Omega_0.$$
Thus, by the definition of $A^{int}$, $\phi\in D(A^{int})$ and $\mu\in \sigma^{int}$.

More precisely, we use the facts that $$\mathrm{Tr}_{ext}(W^D(\Omega))=\mathrm{Tr}_{int}(H^1(\Omega_0))=B^{2,2}_\beta(\Gamma),$$
and thus, the extensions $$E_{ext}: \phi\in B^{2,2}_\beta(\Gamma) \mapsto u\in W^D(\Omega)\hbox{ and } E_{int}: \phi\in B^{2,2}_\beta(\Gamma) \mapsto w\in H^1(\Omega_0)$$
are linear bounded operators. Consequently, $\mu>0$ is an eigenvalue of a Dirichlet-to-Neumann operator with an eigenfunction $\psi\in L_2(\Gamma)$ if and only if  $\mu \phi$ is a normal derivative on $\Gamma$ of $u\in W^D(\Omega)$ or of $w\in H^1(\Omega_0)$, if and only if $\mathrm{Tr}_{ext} u=\phi$ with $\Delta u=0$ weakly on $\Omega$ and $\mathrm{Tr}_{int}w=\phi$ with $\Delta w=0$ weakly on $\Omega_0$, by the uniqueness of the trace and of the normal derivative on $\Gamma$.  Hence, if $\mu\ne 0$, $$\mu\in \sigma^{int}\quad \Longleftrightarrow \quad \mu\in \sigma^{ext}.$$

For $n=2$, $A^{int}$ and $A^{ext}$ are defined in the same way (by $(a_0(\cdot,\cdot),\mathrm{Tr})$), and hence, as in the case $n\ge 3$, the statement
$\sigma^{int}=\sigma^{ext}$ is also a direct corollary of the definitions of the Dirichlet-to-Neumann operators with the continuous extension operators and surjective trace operators mapping to their images.
  
  Now, let us prove Eq.~(\ref{EqT1}). Formula~(\ref{EqT1}) was explicitly proved by Grebenkov~\cite{GREBENKOV-2004} for an annulus p.130. See also~\cite{ARXIV-GIROUARD-2014}. Therefore, it also holds, by a conformal mapping, for domains with regular boundaries. Let us prove it in the general case.   Indeed, since for $n=2$ we have 
  $$0\in\sigma^{int}=\sigma^{ext},\quad \hbox{ and } 0\notin \sigma_S(r).$$
  Moreover, since $\tilde{H}^1(\Omega_{S_r}) \subset H^1(\Omega)$, the functions $u_r\in \tilde{H}^1(\Omega_{S_r})$ can be considered as elements of $H^1(\Omega)$, if outside of $S_r$ we put them  equal to zero. Thus, if $\mu(r)>0$ is an eigenvalue of $A(S)(r)$ in $\Omega_{S_r}$, corresponding to an eigenfunction $\phi\in L_2(\Gamma)$, then 
  for $u_r\in H^1(\Omega)$,  the solution of the Dirichlet Laplacian on $\Omega_{S_r}$, and $u\in H^1(\Omega)$, the solution of the Laplacian with Dirichlet boundary conditions at the infinity (see Remark~\ref{RemLaplExtN2} and Theorem~\ref{PropDtNn2}), we have
  \begin{multline*}
    \forall v \in H^1(\R^n)\quad  \int_{\Gamma}\mu(r) \phi_r \mathrm{Tr}v\dm_d=\int_{\Omega_{S_r}} \nabla u_r\nabla v\dx\\
    \to \int_{\Omega} \nabla u\nabla v\dx =  \int_{\Gamma}\mu \phi \mathrm{Tr}v\dm_d\hbox{ for } r\to +\infty.
  \end{multline*}
  This means that one of the eigenvalues in  the spectrum $\sigma_S(r)$ necessarily  converges towards zero.

 $\Box$

Let us also notice that for the convergence of the series~(\ref{eq/Phi}) on the truncated or the exterior domain, we need to have $\mathds{1}_\Gamma\in \mathcal{D}(A)$. For a Lipschitz boundary $\Gamma$ it was proven in Proposition 5.7 of Ref.~\cite{ARENDT-2015}. In this framework we state more generally 
\begin{proposition}\label{Prop1f}
 Let $\Omega$ be an admissible  exterior domain of $\R^n$ ($n\ge 3$) with a compact $d$-set boundary $\Gamma$, $n-2< d_\Gamma<n$ and let $\Omega_S$ be its admissible truncated domain with $n-2< d_S<n$. Then \begin{equation}\label{DomDefAn3}
      \forall \psi\in L_2(\Gamma)\quad  \exists \phi=A\psi\in L_2(\Gamma),                                                                                                                                  
 \end{equation}
  which also holds for the admissible truncated domains of $\R^2$.

 If $\Omega$ is an admissible exterior domain in $\R^2$ or an admissible domain, bounded by the boundary $\Gamma$ ($n\ge 2$), then 
 \begin{equation}\label{DomDefAn2I}
      \forall \psi \in B^{2,2}_{d/2}(\Gamma)\quad  \exists \phi=A\psi\in L_2(\Gamma).                                                                                                                                 
 \end{equation}
\end{proposition}
\textbf{Proof.}
 Eq.~(\ref{DomDefAn3})  is a corollary of the fact that   the operator $A: L_2(\Gamma)\to L_2(\Gamma)$, considered on $\Omega$ (for $n\ge 3$) and $\Omega_S$  (for $n\ge 2$) respectively, is invertible with a compact inverse operator $A^{-1}$ (since $\lambda=0$ is a regular point by Theorem~\ref{ThSPECTR}).

 For instance, for the exterior case with $n\ge 3$, $\mathds{1}_\Gamma\in L_2(\Gamma)$, thus, for $\lambda=0$, $B_0\mathds{1}_\Gamma \in W^D(\Omega)$, by the well-posedness of the Robin Laplacian exterior problem, and hence, $\mathrm{Tr}(B_0\mathds{1}_\Gamma)=A^{-1}\mathds{1}_\Gamma\in L_2(\Gamma)$.
 
 If $\Omega$ is an exterior domain in $\R^2$ or a bounded domain (see the interior case in Subsection~\ref{SsecDtoNBound}), then for all $u\in H^1(\Omega)$, such that $\Delta u=0$ weakly, there exists unique $\mathrm{Tr} u\in B^{2,2}_\beta(\Gamma)\subset L_2(\Gamma)$ with $\beta=d/2$ (see Theorem~\ref{ThContTrace}), thus for all $\psi \in B^{2,2}_\beta(\Gamma)$ there exists $\phi=A\psi\in L_2(\Gamma)$, as it is stated in Eq.~(\ref{DomDefAn2I}). Consequently, as $1_\Gamma\in  B^{2,2}_\beta(\Gamma)$, we have $\mathds{1}_\Gamma\in \mathcal{D}(A)$.
$\Box$

%
\section*{Acknowledgments} We would like to thank Claude Bardos for useful discussions about the subject and  Denis Grebenkov for pointing the physical meaning of the problem.


          \label{bib:sec}


%

\end{document}